\documentclass[12pt]{article}

\usepackage{epsfig}

\usepackage{amssymb}

\newtheorem{theorem}{Theorem}[section]
\newtheorem{lemma}[theorem]{Lemma}
\newtheorem{corollary}[theorem]{Corollary}
\newtheorem{remark}[theorem]{Remark}

\newtheorem{definition}[theorem]{Definition}
 \newenvironment{proof}{{\it Proof:\/}}{$\Box$\vskip 0.08in} 

\newtheorem{problem}[theorem]{Problem}
\newtheorem{example}[theorem]{Example}

\newtheorem{formulla}[theorem]{}

\begin{document}
\begin{center} {\large\bf
Homotopy and q-homotopy skein modules of 3-manifolds: an example
in Algebra Situs.} \end{center}

\begin{center}
J\'ozef H.~Przytycki\footnote{Supported by NSF-DMS-98089555.} \end{center}

\begin{quotation}
\begin{center} Dedicated to my teacher Joan Birman on her 70'th birthday.
\end{center} \end{quotation}

\vspace{0.1in}\ \\ 
\begin{abstract}
{\it Algebra Situs} is a branch of mathematics which has its roots
in Jones' construction of his polynomial invariant of links and
Drinfeld's work on quantum groups. It encompasses the theory of
quantum invariants of knots and 3-manifolds, algebraic topology
based on knots, operads, planar algebras, $q$-deformations,
quantum groups, and overlaps with algebraic geometry,
non-commutative geometry and statistical mechanics.

{\it Algebraic topology based on knots} may be characterized as a
study of properties of manifolds by considering links
(submanifolds) in a manifold and their algebraic structure. The
main objects of the discipline are {\it skein modules}, which are
quotients of free modules over ambient isotopy classes of links in
a manifold by properly chosen local (skein) relations.

We concentrate, in this lecture, on one relatively simple example
of a skein module of 3-manifolds -- the $q$-homotopy skein module.
This skein module already has many ingredients of the theory:
algebra structure, associated Lie algebra, quantization, state
models...
\end{abstract}

\section{Introduction}\label{1}

Algebra Situs\footnote{This part of the paper is based on the talk
{\it Algebraic topology based on knots: a case study in the
history of ideas} given at a Conference in Low-Dimensional
Topology in Honor of Joan Birman's 70th Birthday; Columbia
University, March 14--15, 1998.}
 is a branch of mathematics which has its roots in
Jones' construction of his polynomial invariant of links, Jones
polynomial, and Drinfeld's work on quantum groups. It encompasses
theory of quantum invariants of knots and 3-manifolds, algebraic
topology based on knots, operads, $q$-deformations, quantum
groups, and overlaps with algebraic geometry, non-commutative
geometry and statistical mechanics.

Algebraic topology based on knots may be characterized as a study
of the properties of manifolds by considering the space of links
(submanifolds) in a manifold and its algebraic structure. The main
objects of the discipline are skein modules, which are quotients
of free modules over ambient isotopy classes of links in a
manifold by properly chosen local (skein) relations. Of course,
this is not a complete definition of the field, which has its
purely algebraic component (skein algebras of groups), higher
manifold generalization and rich internal structure, but at least
it gives the idea of our subject.

In searching for a starting point of the theory one should
consider Listing book (1847), Dedekind and Weber's paper (1882),
and Poincar\'e's paper ``Analysis Situs" (1895). In knot theory,
skein modules (building blocks of algebraic topology based on
knots) have their origin in the observation by Alexander (1928)
that his polynomials of three links $L_+, L_-$ and $L_0$ in $S^3$
(see Fig.2.1) 
are linearly related. This line of research was continued by
Conway (linear skein, 1969). In graph theory the idea of forming a
ring of graphs and dividing it by an ideal generated by local
relations was developed by Tutte in his 1946 PhD thesis. The
theory of Hecke algebras, as introduced by Iwahori (1964), is
closely connected to the theory of skein modules. Another
connection can be found in the Temperley-Lieb algebra (1971). The
main motivation for skein modules was the discovery/construction
of the Jones polynomial (1984). Skein algebras of groups use rich
ideas of Poincar\'e (1884) , Vogt (1889), Fricke and Klein (1897)
and the school of Magnus (e.g. ``Rings of Fricke characters",
1980).

Joan Birman introduced me to the world of knots and braids, to the
work of her advisor W.Magnus, and her grand-advisor
M.Dehn.\footnote {See Fig.1.1 for genealogical table.} Her work
is continued by her students and grand-students (one of the best
recent results related to algebraic topology based on knots was
obtained by A.Sikora \cite{Si-2}).

We concentrate, here, on one relatively simple example of a skein
module of 3-manifolds -- the q-homotopy skein module. This skein
module already has many ingredients of the rich theory: algebra
structure, associated Lie algebra, Hopf algebra, quantization,
state models, relation to graph theory...

\ \\
\centerline{\psfig{figure=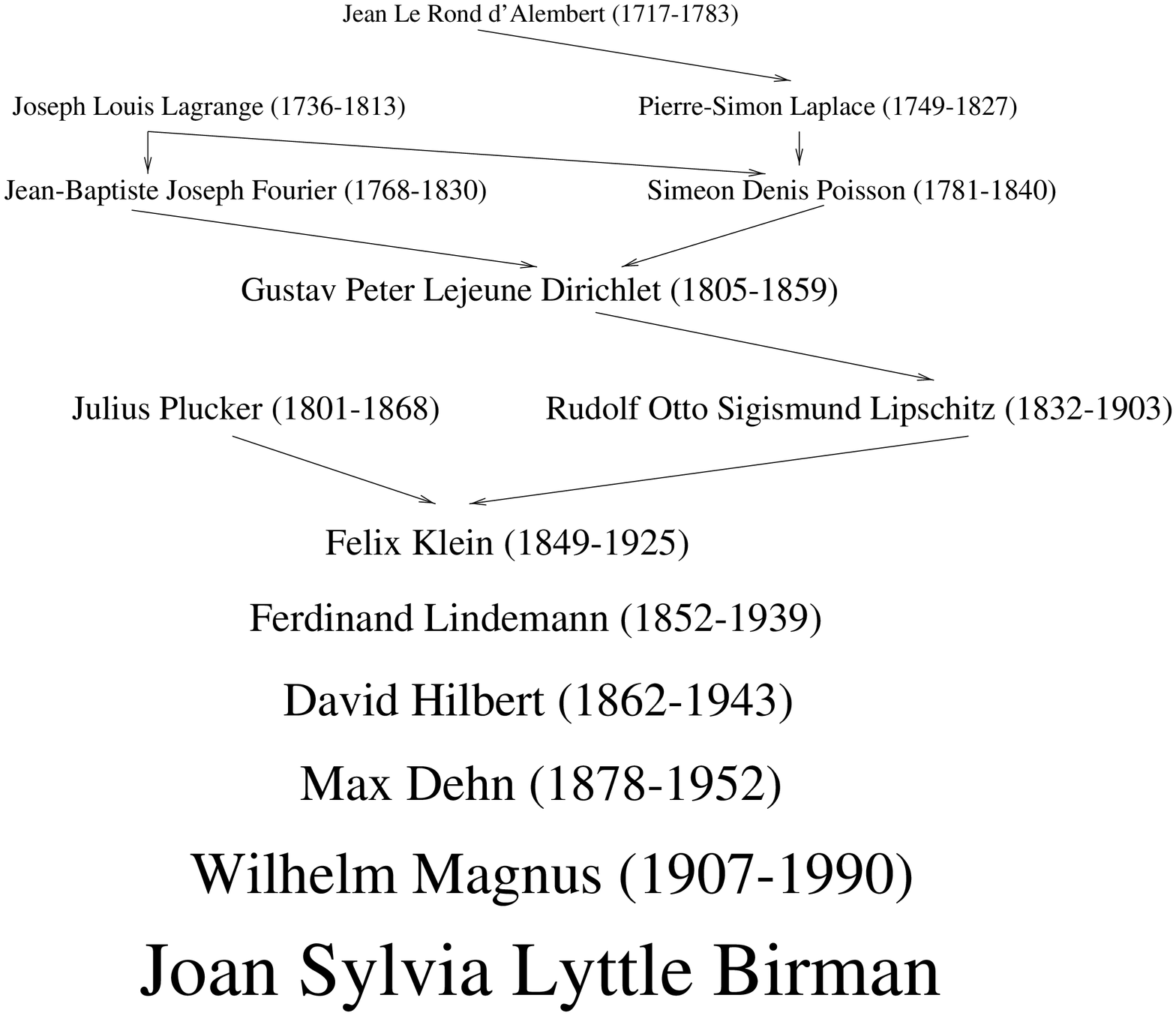,height=13.5cm}}
\ \\
\centerline{Students of J.Birman:}
Tara Brendle, Nathan Broaddus, Abhijit Champanerkar, Zung He Chen,
Richard Fein, Liz Finkelstein, Tat San Fun, Matt Greenwood, Hessam
Hamidi-Tehrani, Efstratia Kalfagianni, Marcello Kupferwasser, John
D. McCarthy, Ka Yi Ng, Radu Popescu, Jerry Powell, J\'ozef H.
Przytycki, Ted Stanford, Rolland J. Trapp, Nancy Wrinkle, Peijun
Xu, Matt Zinno,..., Keiko Kawamuro.
\begin{center}
Figure 1.1.\ \ Teachers and students of Joan Birman. \end{center}
\newpage
\section{Definitions}
 One of the earliest skein modules considered\footnote{Sections 2-4
are based on notes
written in October and November of 1991 when I was participating
in the Topology Semester at UTK \cite{Pr-3}. I would like to thank
Morwen Thistlethwaite and the University of Tennessee for
hospitality.} was
the homotopy skein module, where self-crossings were ignored and
mixed crossings were resolved using the Alexander-Conway skein
relation \cite{H-P-1}. In the fall of 1991 the $q$-analogue
version of the homotopy skein module was first considered
\cite{Pr-3} and showed to distinguish some links with the same
Jones-Conway (Homflypt) polynomial. In the case of $M=F\times I$
we gave in \cite{H-P-1} the precise description of the homotopy 
skein module and
showed that as an algebra it is isomorphic to the universal
enveloping algebra of the Goldman-Wolpert Lie algebra of curves on
the surface $F$, \cite{Gol}. For the $q$-analogue version of the
homotopy skein module we prove that for $M=F\times I$ and
$\pi_1(F)$ abelian, the skein module is free. We define a
$q$-version of the Goldman-Wolpert Lie algebra. We also give the
formula for an element of the $q$ homotopy skein module
represented by a link in $S^3$ in terms of linking numbers of the
components of the link. We show that the module (algebra), for
$M=S^3$, is equivalent to the dichromatic polynomial of the graph
associated to the link. We generalize these to surfaces $F$ with
an abelian fundamental group. We show that for $\pi_1(F)$
nonabelian the q-homotopy skein module has torsion.

The $q$-homotopy skein module ${\cal H \cal S}^q(M)$ is a
$Z[z,q^{\pm 1}]$-module associated to an oriented 3-manifold, $M$,
according to the general scheme described in \cite{Pr-2,H-P-2}. It
generalizes the homotopy skein module ${\cal H \cal S}(M)$
considered in \cite{H-P-1} and it can be thought of as a
$q$-analogue of ${\cal H \cal S}(M)$. We will follow \cite{H-P-1}
closely in our description. Let $M$ be an oriented 3-manifold. The
$q$-homotopy skein module ${\cal H \cal S}^q(M)$ is defined as
follows. \begin{definition} Let ${\cal L}^h$ be the set of all
oriented links in $M$ modulo link homotopy (i.e. we ignore
self-crossings), including the empty link , $T_0$. Let $R=Z[q^{\pm
1},z]$ and $R{\cal L}^h$ be the free $R$ module generated by
${\cal L}^h$. Let $S$ be a submodule of $R{\cal L}^h$ generated by
the homotopy skein expressions\ \ $q^{-1}L_+ -qL_- -zL_0$, where
$L_+, L_-$ and $L_0$ are three oriented links in $M$, which are
identical except inside a small ball where they appear instead as
shown in Figure 2.1.\\ \ \\ \ \\
\centerline{\psfig{figure=L+L-L0.eps,height=2.2cm}}
\begin{center}
Figure 2.1
\end{center}
\ \\ Additionally we assume that the two strings of $L_+$ (or
$L_-$) involved in the crossing in Fig.2.1 belong to different
components of $L_+$ (or $L_-$), that is, we deal with a mixed
crossing. We define the $q$-homotopy skein module to be the
quotient ${\cal H \cal S}^q(M)= R{\cal L}^h/S$.

If we do not allow the empty link, we get the reduced $q$-homotopy
skein module ${\cal H} \tilde{S}^q(M)$ (and we have ${\cal H \cal
S}^q(M)= {\cal H} \tilde{S}^q(M)\oplus R$).
\end{definition}

\begin{remark}\label{2.2}
It may also be convenient to define ${\cal H \cal S}^q(M)$,
equivalently, starting from the set of oriented links, $\cal L$,
in $M$ and quotient $R\cal L$ by the submodule generated by skein
relations $L_+ - L_-$ in the case of a self-crossing and
$q^{-1}L_+ -qL_- -zL_0$ in the case of a mixed crossing.
\end{remark}
The $q$-homotopy skein module shares with other skein modules
several useful elementary properties, like the Universal
Coefficients property, several functorial properties, etc. We will
discuss these in Section 5. In the next two sections we analyze
the $q$-homotopy skein module of classical links (that is $M=S^3$)
and show that it can be interpreted as a polynomial, which we call
the {\it homotopy polynomial}. The homotopy polynomial corresponds
to (a variant of) the dichromatic polynomial of an associated graph,
and depends only on the linking numbers between components of the
link.

\section{${\cal H \cal S}^q(S^3)$ and a homotopy polynomial of classical 
links}

It is relatively easy to show that ${\cal H \cal S}^q(S^3)$ is
freely generated by unlinks $T_0,T_1, T_2, T_3,...$, where $T_i$
denotes the unlink of $i$ components, and that for a given
$n$-component link $L\subset S^3$, its presentation as
$L=w_0(q)T_n + w_1(q)zT_{n-1} +...+w_{n-1}(q)z^{n-1}T_{1}$ depends
only on the linking numbers between components of $L$ (compare
\cite{H-P-1}; Theorem 1.2). In Sections 6 we give detailed proof
of much more general result. However, the case of $M=S^3$ is of
special interest as ${\cal H \cal S}^q(S^3)$ is equivalent to a
classical object -- the dichromatic polynomial, properties of
which are very well understood. More precisely, we can put
$T_i=t^i$ and interpret $L \in {\cal H \cal S}^q(S^3)$ as
$L=HP_L(q,t,z)T_0$, where $HP_L(q,t,z) \in Z[q^{\pm 1},t,z]$; in
the previous notation $HP_L(q,t,z)= w_0(q)t^n + w_1(q)zt^{n-1}
+...+w_{n-1}(q)z^{n-1}t$. We call the polynomial $HP_L(q,t,z)$ the
{\it homotopy polynomial} of a link $L$. It can be interpreted as
a dichromatic polynomial of the weighted graph associated to $L$,
and the set of linking numbers between components of $L$ can be
recovered from the coefficient $w_1(q)$.\footnote{In \cite{H-P-1}
we didn't write the closed formulas for $w_k(1)$. I noticed
Formula 3.4 shortly after \cite{H-P-1} was published and
generalized to $w_k(q)$ \cite{Pr-3}. Formula for $w_k(1)$ was also
independently discovered by A.Sikora in his Master Degree Thesis
\cite{Si-1}.}

To formulate the main result of this section we need some preliminary
definitions. For a link $L= K_1 \cup ... \cup K_n$
in $S^3$ we denote by $[{l}_{i,j}]$ its linking matrix where
${l}_{i,j}=lk(K_i,K_j)$. Let $E$ be the set of all pairs 
$(i,j)$, $i\neq j$, $1\leq i,j \leq n$. We consider the notion of 
a cycle in $E$. 
The meaning of a cycle can be best
explained by considering the complete graph $F_n$ of $n$ vertices,
$1,2,3,...,n$, and edges $(i,j)$ ($i\neq j$). Cycles (i.e. simple
closed edge-paths) in $F_n$ determine cycles in $E$. For $S\subset E$
where $S$ does not contain a cycle,
let $A_S$ denote the subset of $E-S$ such that
$(i,j)\in A_S$ if and only if either $S\cup \{(i,j)\}$ has no
cycle or, otherwise, if $C$ is the unique cycle in $S\cup
\{(i,j)\}$ containing $(i,j)$ then $(i,j)$ is not the first
element of $C$ with respect to lexicographical order of pairs
$(i,j)$. \footnote{ $A_S$ correspond to the set of externally
inactive elements in the sense of Tutte; compare \cite{Tut-4,Tra}.} $|S|$
denotes the cardinality of $S$.

\begin{theorem}\label{3.1}
$T_0,T_1, T_2, T_3,...$ form a free basis of ${\cal H \cal
S}^q(S^3)$ and for $L= K_1 \cup ... \cup K_n$ 
one has the formula:
\end{theorem}
{\begin{formulla}\label{3.2} $$L=\sum_{S} q^{\Sigma_{(i,j)\in A_S}
2{l}_{i,j}}z^{|S|}T_{n-|S|} \prod_{(i,j)\in
S}\frac{q^{2{l}_{i,j}}-1}{q-q^{-1}} $$
where the sum is taken over all subsets $S$ of $E$ which do
not contain a cycle.

Equivalently we can write:
$$ HP_L(q,t,z)=\sum_{S} q^{\Sigma_{(i,j)\in A_S}
2{l}_{i,j}}z^{|S|}t^{n-|S|} \prod_{(i,j)\in
S}\frac{q^{2{l}_{i,j}}-1}{q-q^{-1}} $$ 
where as before the sum is taken over all subsets $S$ of $E$ which do
not contain a cycle.
\end{formulla}

\begin{corollary}\label{3.3}
If $q=1$, then the formula in 3.1 reduces to:
\end{corollary}
{\begin{formulla}\label{3.4}
$$L=\sum_S z^{|S|}T_{n-|S|} \prod_{(i,j)\in S}{l}_{i,j}$$ \end{formulla}} 

{\it In particular $w_k(1)$ is equal to $\sum_{S_k}\prod_{(i,j)\in
S_k}{l}_{i,j}$, where $S_k$ is the set of all $k$-element subsets
of $E$ which do not contain a cycle. Compare \cite{H-P-1},
formulas of Part 1.}

The proof of Theorem 3.1 is not very difficult, but we can omit it
totally\footnote {In Section 6 we prove a generalization of the
first part of Theorem 3.1.} by showing that Formula 3.2 can 
be interpreted as a formula for dichromatic polynomial of
the signed (or weighted) graph associated to $L$.

We will consider three, closely related, graphs: $G(L), G_1(L)$
and $G_2(L)$ and their dichromatic polynomials. In all of them
vertices correspond to components of $L$.

Consider the following signed graph $G(L)$:\ \ vertices of $G(L)$
correspond to components of $L$ and vertices $v_i,v_j$ ($i\neq j$)
are joined by $|lk(L_i,L_j)|$ edges of sign equal to
$sign(lk(L_i,L_j))$.

Let $R(G)$ be a dichromatic polynomial of a given signed graph $G$
defined recursively by the rules:
\begin{enumerate}
\item[(i)] $R(\underbrace{\bullet \bullet ...\bullet}_n)
=t^n$,
\item[(ii)] if $e_{\pm}$ is not a loop then:
$$ R(G)=q^2R(G-e_+) + qzR(G/e_{+})$$ $$ R(G)=q^{-2}R(G-e_-) -
q^{-1}zR(G/e_{-}),$$
\item[(iii)] If $e$ is a loop then $R(G)=R(G-e)$.
\end{enumerate}
Remark.\\ Notice that if $e_+$ and $e_-$ are two edges joining the
same endpoints, then $R(G-e_+ - e_-) = R(G)$ because, either $e_+$
and $e_-$ are loops and the equality follows from (iii), or we
have:\\ $R(G)= q^2R(G-e_+) + qzR(G/e_{+})= q^2(q^{-2}R(G-e_+ -
e_-)- q^{-1}zR((G-e_+)/e_-))+ qzR(G/e_{+}) = R(G-e_+ - e_-)$ as
$e_-$ is a loop in $G/e_+$ so $R(G/e_+)=R((G-e_+)/e_-)$.

It is a standard fact that $R(G)$ is well defined, in particular
one has the state model formula (compare \cite{F-K,Tra,P-P}). One
can immediately see a validity of Formula 3.5 when 
one notices that (iii) can be rewritten as:\\ 
$$ R(G)= \left \{ \begin{array}{ll}
q^2R(G-e_+) + qz\frac{q^{-1}-q}{z}R(G/e_{+}) & \\ q^{-2}R(G-e_-) -
q^{-1}z\frac{q^{-1}-q}{z}R(G/e_{-}) & \end{array} \right. $$

\begin{formulla}\label{3.5}
$$R(G)=\sum_{S\in E(G)}t^{p_0(G:S)}(\frac{q^{-1}-q}{z})^{p_1(G:S)}
(-1)^{|S|_-}z^{|S|}q^{2(|E-S|_+ - |E-S|_-)+|S|_+-|S|_-}$$
\end{formulla}
where $(G{:}S)$ is the subgraph of $G$ which includes all vertices
of $G$ but only edges of $S$. $p_0(G)$ is the number of components
of $G$, and $p_1(G)$ is the cyclomatic number of $G$ (i.e. the
first Betti number). $V=V(G)$ denotes the set of vertices of $G$
and $|V|$ the cardinality of $V$. $E=E(G)$ denotes the set of
edges of $G$ and $|E|$ (resp. $|E|_+$ or $|E|_-$) denotes the
cardinality of $E$ (resp. the number of positive or negative edges
in $E$). In particular $p_1(G)= |E| - |V| + p_0(G)$.

\begin{remark}
\begin{enumerate}
\item[(a)]
If $Q(G;t,z')$ is Traldi's version of the dichromatic polynomial
(\cite{Tra}), then $$R(G)= q^{2(|E|_+ -|E|_-)}Q(G;t,z'),$$ for
$z'=\frac{q^{-1}-q}{z}$ and Traldi's weight $w(e)$ of $e\in E$ is
defined by $w(e_+)=q^{-1}z$ and $w(e_-)=-qz$.
\item[(b)] If $<G>_{\mu,A,B}$ is the Kauffman bracket of $G$ then
(compare \cite{P-P}, Lemma 5.2):
$$<G>_{\mu,A,B}=\mu^{-1}(-A^3)^{|E|_+ -|E|_-}R(G)_{z=-i}$$ for
$\mu =t$, $A= (iq)^{-\frac{1}{2}}$, $B= (iq)^{\frac{1}{2}}$ and if
$R(G)_{z=-i}=\sum b_j(q)t^j$ then $R(G)=\sum
b_j(q)t^j(iz)^{|V|-j}$.
\end{enumerate}
\end{remark}

\begin{theorem}\label{3.7} Let $\hat w : {\cal L} \to {\cal G}$ be a map from 
the set of links in $S^3$ to the set of signed graphs, given by 
$\hat w(L) =G(L)$. Then $\hat w $ yields an algebra isomorphism 
$w: {\cal HS}^q(S^3) \to Z[q^{\pm 1},z,t]$, where $w(T_i)=t^i$ and
$w(L)=R(G(L))$. The product in ${\cal HS}^q(S^3)$ is given by the
disjoint sum of links. Furthermore $L=R(G(L))T_0$ in ${\cal
HS}^q(S^3)$.
\end{theorem}
\begin{proof}
Because $T_i$'s clearly generate ${\cal HS}^q(S^3)$, it suffices
to compare defining properties (i)-(iii) of $R(G)$ with the
definition of ${\cal HS}^q(S^3)$. The condition (ii) corresponds
to the homotopy skein relation.
\end{proof}
Notice that we have proven the first part of Theorem 3.1. Formula
3.2 can be deducted from 3.5, but we will not present it here.
Instead we will consider a graph $G_1(L)$ yielded by $L$ and from
its dichromatic polynomial we will derive 3.2.

We can construct a graph $G_1$ out of any signed graph $G$.
\begin{definition}
Let $G$ be a signed graph, then $G_1$ is a weighted graph (i.e. a
graph with a function $f: E(G_1) \to Z$), with no loops and no
multi-edges, obtained from $G$ by deleting its loops and
replacing a multi-edge by a single edge with the weight equal to
the sum of signs of edged in the multi-edge. Deleting, $G_1-e$,
and contracting, $G_1/e$, are operations on $G_1$ with the usual
meaning with the convention that whenever a multi-edge is
created in $G_1/e$ then it is replaced by a single edge with a
weight being the sum of weights of components of the multi-edge.

Now we define the dichromatic polynomial $R_1(G_1)$ by the rules:
\begin{enumerate}
\item[(i)] $R_1(\underbrace{\bullet \bullet ...\bullet}_n)
=t^n$,
\item[(ii)] $ R_1(G_1)=q^{2f(e)}R_1(G_1-e) + \frac{q^{2f(e)}-1}
{q-q^{-1}}zR_1(G_1/e)$, where $f(e)$ is the weight of $e$. Notice
that $e$ is never a loop (as we delete edges of weight zero).\
Observe also that for $q=1$, (ii) reduces to:
\item[(ii')] $ R_1(G_1)_{q=1} = R_1(G_1-e)_{q=1} + f(e)zR_1(G_1/e)_{q=1}$.
\end{enumerate}
\end{definition}
It is easy to check that if $e$ is an edge of $G_1$ with weight
$f(e)$ and $G'$ is obtained from $G_1$ by changing the weight of
$e$ to $f(e)-1$ then $R_1(G_1)=q^2R_1(G') + qzR_1(G_1/e)$. From
these we get:
\begin{lemma}
$R(G)=R_1(G_1)$.
\end{lemma}
The following theorem of Tutte \cite{Tut-4} (compare \cite{Tra} or
\cite{Za}) yields Theorem 3.1.
\begin{theorem}\label {3.10}
Let $G_1$ be a connected weighted graph. Given an arbitrary linear
ordering of edges in $E(G_1)$ and a spanning tree $T$ of $G_1$, an
edge $e$ in $T$ is called internally active with respect to $T$ if
it precede all other edges of $G$ whose end vertices lie in
different components of $T-e$. An edge $e$ not in $T$ is called
externally active with respect to $T$ if it precedes all other
edges of $T$ that lie in the unique cycle determined by $T$ and
$e$. Then:
\end{theorem}
\begin{formulla}\label{3.11}
$$R_1(G_1)= t\sum_{T}\prod_{e\in EI}q^{2f(e)} \prod_{e\in
II}(z\frac{q^{2f(e)}-1}{q-q^{-1}}) \prod_{e\in IA}(tq^{2f(e)}+
z\frac{q^{2f(e)}-1}{q-q^{-1}})$$ where the sum is taken over the
set of all spanning trees $T$ of $G_1$ and the three products are
taken, respectively, over the set of edges which are externally
inactive (EI), internally inactive (II) and internally active (IA)
with respect to $T$.
\end{formulla}
Let $|S|$ denote the number of edges of $S \subset E(G_1)$, and
$n$ number of vertices of $G_1$. If we multiply out the last
product of Formula 3.11 (having in mind the defining relations
(i)-(iii) for $R_1(G_1)$) we get\footnote{ We can visualize this
by applying recursive relations to the edges of the trees.}:
\begin{formulla}\label{3.12}
$$R_1(G_1)=\sum_{S}t^{n-|S|}z^{|S|}\prod_{e\in A_S}q^{2f(e)}
\prod_{e\in S}(\frac{q^{2f(e)}-1}{q-q^{-1}})$$ \end{formulla}
where the sum is taken over all forests $S$ of $G_1$ (that is,
subgraphs without cycles), $A_S$ denotes the set of externally
inactive edges with respect to $S$ (i.e. $e$ in $E-S$ is inactive
if $S\cup e$ is a forest or $S\cup e$ contains a (unique) cycle
and $e$ does not precedes all other edges of $C$). Formula 3.2 now
follows from 3.12 if edges of $G_1$ are lexicographically ordered.
The proof of Theorem 3.1 is complete.

We will end this section by introducing another related graph,
$G_2$, and showing that we can interpret Jones-Conway (Homflypt)
polynomial of links in $S^3$ as a dichromatic polynomial of a
related graph.
\begin{definition}\label{3.13}
Let $G_2$ be a signed graph obtained from $G$ by doubling each (signed) edge 
of $G$; that is, {\psfig{figure=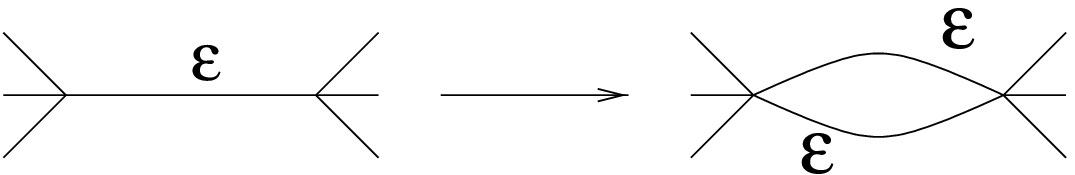,height=0.5cm}} . Let 
$R_2(G_2)$ be a polynomial of $G_2$ given by the rules:
\begin{enumerate}
\item[(i)] $R_2(\underbrace{\bullet \bullet ...\bullet}_n) =t^n$,
\item[(ii)] if $e$ is not a loop and $\epsilon (e)$ denote a sign of the edge 
$e$,
then: $$ R_2(G)=q^{2\epsilon (e)}R_2(G_2(\epsilon (e))) + \epsilon
(e) q^{\epsilon (e)}zR_2(G_2/e)$$ where $G_2(\epsilon (e))$
denotes the graph obtained from $G_2$ by changing the sign of the
edge $e$.
\item[(iii)] If $e$ is a loop then $R_2(G_2)=R_2(G_2-e)$.
\end{enumerate}
\end{definition}
We do not claim that $R_2(\ )$ is defined for every signed graph,
but it is defined for $G_2$ constructed from $G$, as above. In
fact we have:
\begin{lemma}\label{3.14}
$R(G)=R_2(G_2)$.
\end{lemma}
For any plane signed graph $G$, we can associate the link diagram
$D(G)$ (matched diagram) as shown in Fig. 3.1 below\footnote{ We
can also say that $D(G)$ is obtained from $G_2$ by the standard
(from P.G.Tait times) construction -- median graph diagram, with
the properly chosen orientation and crossing resolution.}: \ \\ \
\\
\centerline{\psfig{figure=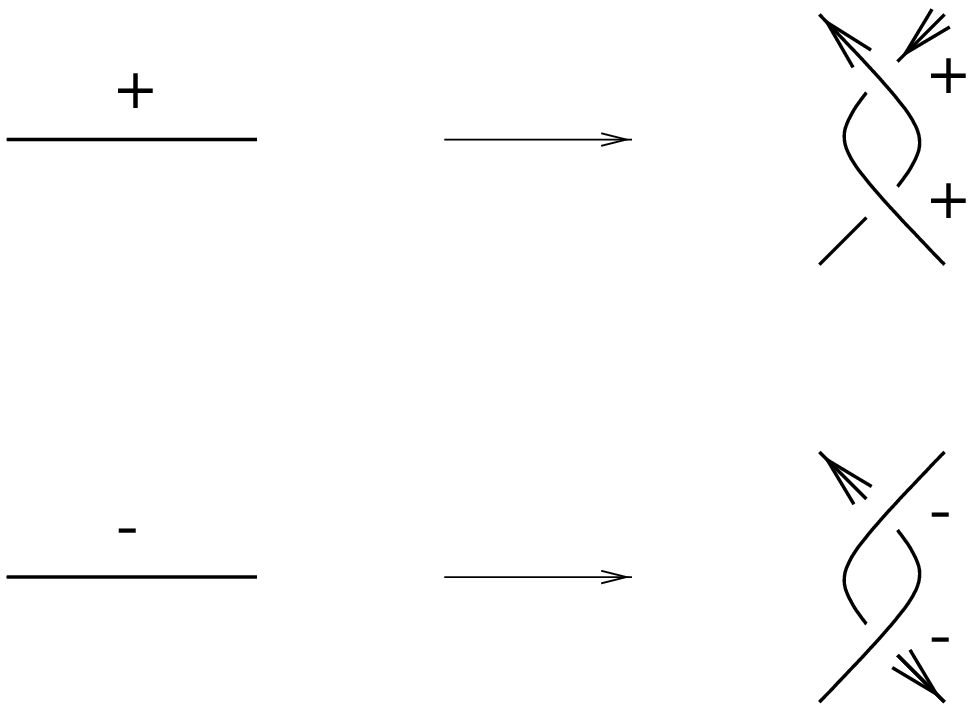,height=3.0cm}}
\begin{center}
Figure 3.1
\end{center}
\ \\ Then the dichromatic polynomial of $G$ and the Jones-Conway
(Homflypt) polynomial of $D(G)$, $P_{D(G)}(v,z)$, are related as
follows:
\begin{lemma}\label{3.15}
$P_{D(G)}(v,z)= t^{-1}R(G)$ for $v=q$ and $t=\frac{v^{-1}-v}{z}$.
Thus we can recover $P_{D(G)}(v,z)$ out of $R(G)$. Conversely if
$P_{D(G)}(v,z)=\sum_{i=0}^m a_{-m +2i}(v)z^{-m+2i}$, where
$m=com(D(G))-1= |V(G)|-1$, then $R(G)=t\sum_{i=0}^m \frac{a_{-m
+2i}(q)}{(q^{-1}-q)^{m-i}}t^{m-i}z^i$. Here $com(D(G))$ denotes
the number of components of the link diagram $D(G)$.
\end{lemma}
\section{Examples}
As noted in Section 3, the $q$-homotopy skein module, ${\cal
HS}^q(S^3)$, can be identified with the ring $Z[q^{\pm 1},t,z]$.
Thus the class of a link $L$ in ${\cal HS}^q(S^3)$ defines the
homotopy polynomial, $HP_L(q,t,z)\in Z[q^{\pm 1},t,z]$, which
depends exclusively on the linking matrix $[{\ell}_{i,j}]$ of $L$.
It is worth comparing $HP_L(q,t,z)$ with the Jones-Conway
(Homflypt) polynomial $P_L(v,z)$. It is well known (compare
\cite{L-M,Pr-1,Si-1}) that $HP_L(1,t,z)$ (i.e. $q=1$) is
determined by $P_L(v,z)$. This is not true, however, for 
the more general polynomial $HP_L(q,t,z)$. 
In the following example, we use links described by J.Birman
\cite{Bir}.

\begin{example}\label{4.1}
The 3-component links from Fig. 4.1 share the same Jones-Conway
polynomial \cite{Bir} but they have different homotopy
polynomials.\\ Namely, $HP_{L_1}(q,t,z) = q^{6}t^3 +
(q^{-1}+q+q^3+q^5-q^7)zt^2 -(1+q^2+q^4+q^6)z^2t $ and
$HP_{L_2}(q,t,z) = q^{6}t^3 + (q+2q^3+2q^5-q^7-q^9)zt^2
-(q^4+2q^6+q^8)z^2t $.
\end{example}
\ \\
\centerline{\psfig{figure=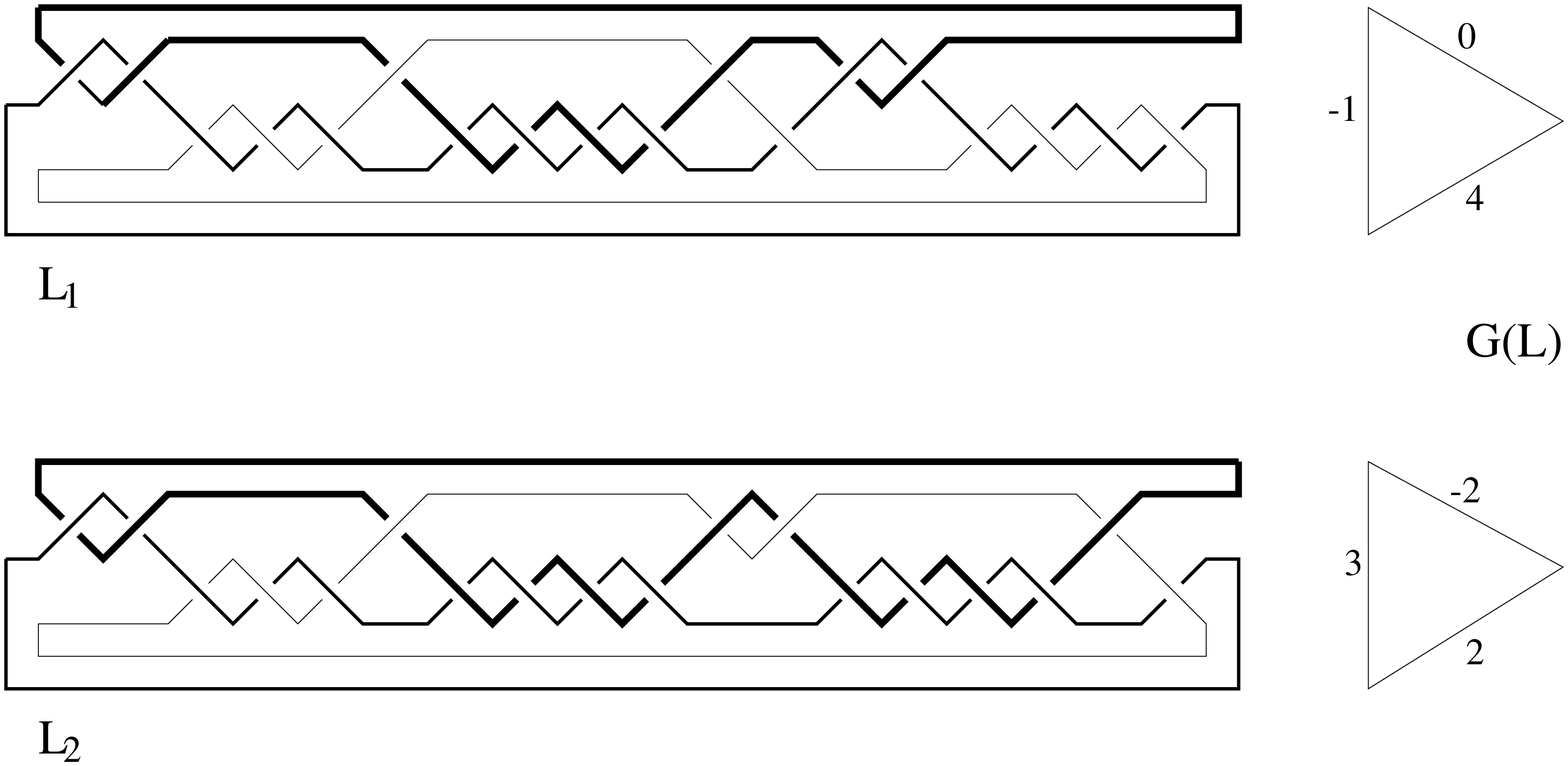,height=4.7cm}}
\begin{center}
Fig. 4.1.
\end{center} Another interesting example is also due to Joan Birman, this time
it deals with two different closed three braids (one being the
mirror image of the other).
\begin{example}\label{4.2}
The 3-component links from Fig. 4.2 share the same Jones-Conway
polynomial \cite{Bir} but they have different homotopy
polynomials. Namely, the coefficient, $w_1(q)$, of $zt^2$ in
$HP_{L}(q,t,z)$ is equal to $-q^3-q+2q^{-1}$, and in $HP_{\bar
L}(q,t,z)$ is equal to $q^{-3} +q^{-1} - 2q$. \footnote{ Generally
for a link $L$ and its mirror image $\bar L$ one has: $HP_{\bar
L}(q,t,z) = HP_{L}(-q^{-1},t,z)$.}
\end{example}
\ \\
\centerline{\psfig{figure=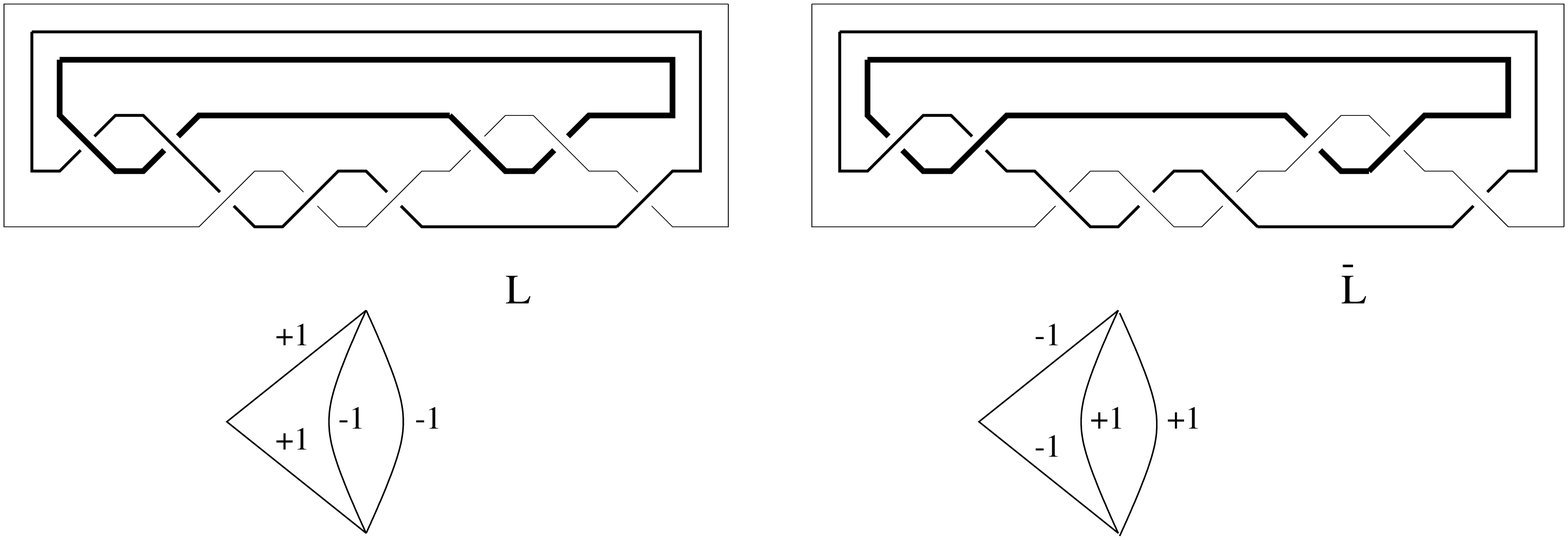,height=4.8cm}}
\begin{center}
Fig. 4.2. \end{center} The polynomial $HP_{L}(q,t,z)$ allows us to
recover all linking numbers of $L$ (as a set with multiplicities);
see Theorem 4.3, but is not sufficient to recover the linking
matrix. The simplest example is shown in Fig. 4.3.\ \\
\centerline{\psfig{figure=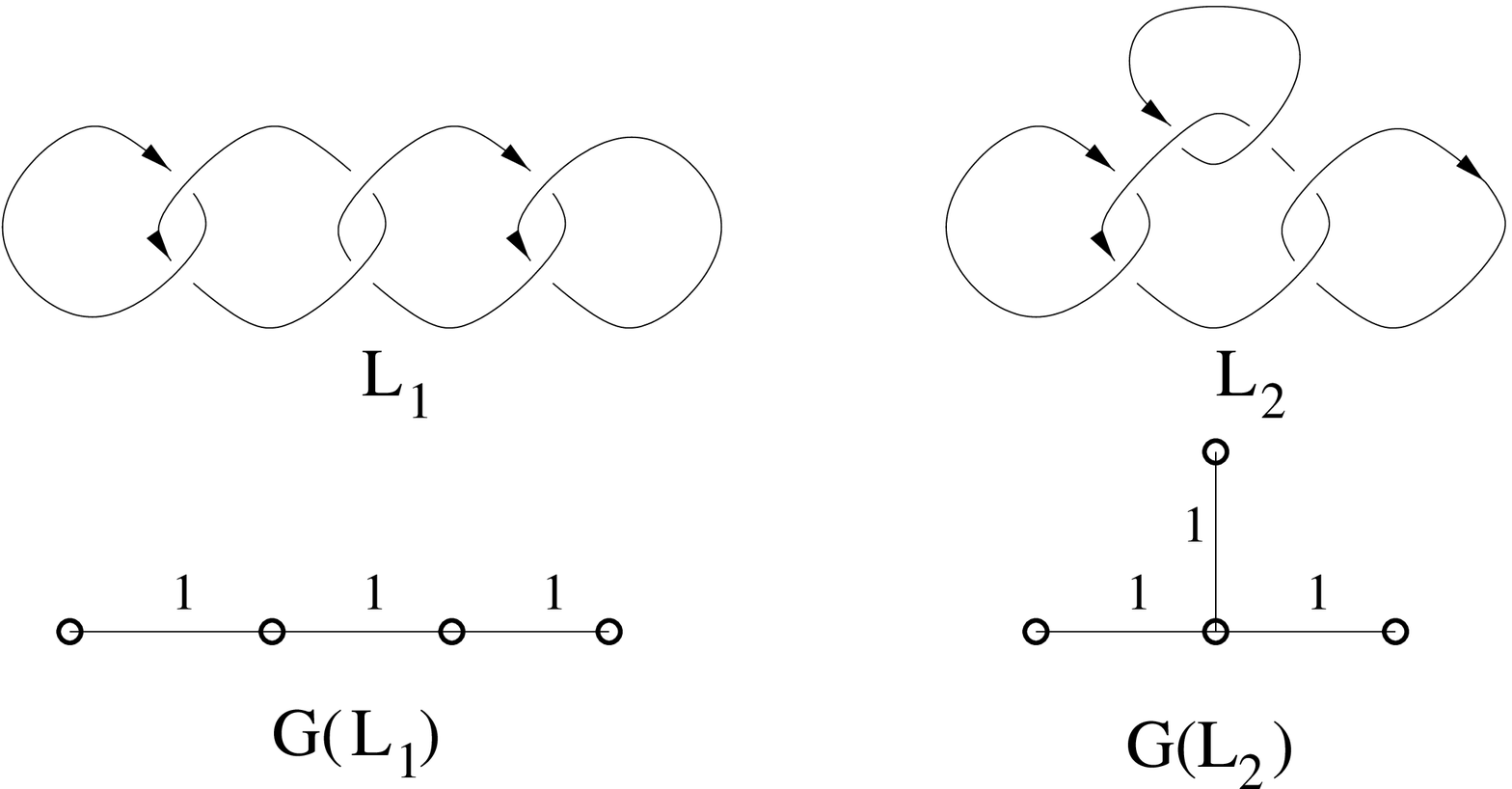,height=4.5cm}}
\begin{center}
Fig. 4.3.
\end{center}
Here $HP_{L_1}(q,t,z) = HP_{L_2}(q,t,z) = t(q^2t+qz)^3$.

This example is, in a sense, trivial because the graphs $G(L_1)$
and $G(L_2)$ are 2-isomorphic. The first nontrivial example is
based on the example of M.C.Gray obtained around 1933 (see
\cite{Tut-2}). The graphs $G(L_1)$ and $G(L_2)$ of Fig. 4.4 are
that of Gray. It is an open question how much of the graph can be
recovered from its dichromatic polynomial and whether there are
some ``elementary moves" on graphs which link different graphs
with the same dichromatic polynomial. Most of examples known today
are based on the ``rotors" idea of Tutte
(\cite{BSST,Tut-2,Tut-3,APR,T,Jo,Pr-5}).

\ \\
\centerline{\psfig{figure=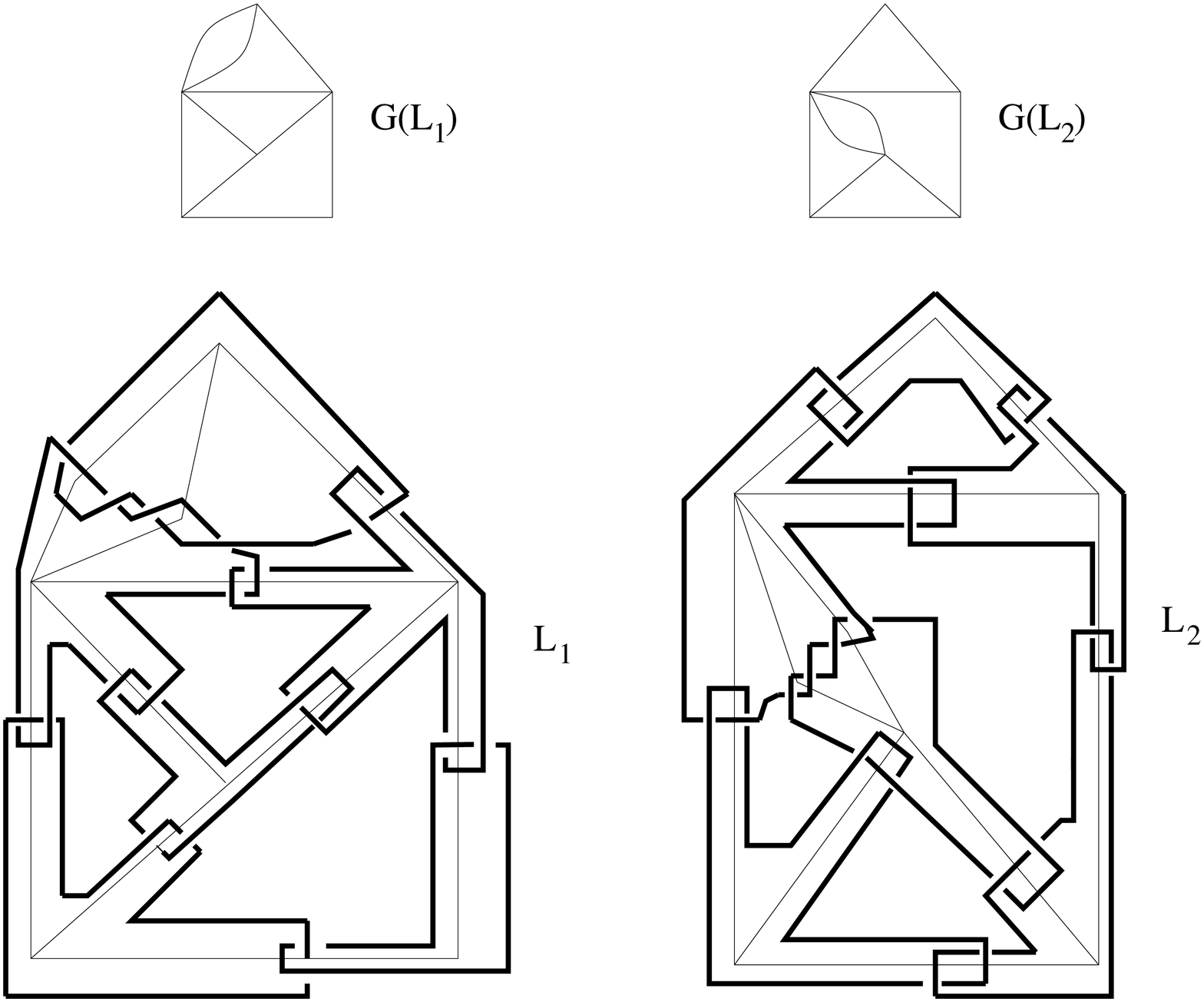,height=8.3cm}}
\begin{center}
Fig. 4.4.
\end{center}
We will now show that $HP_{L}(q,t,z)$ allows us to recover the set
of linking numbers $\{{\ell}_{i,j}\}$.
\begin{theorem}\label{4.3}
Let $L$ be an oriented link on $n$ components in $S^3$ and
$$HP_{L}(q,t,z) = \sum_{i=0}^{n-1}w_i(q)z^it^{n-i}\ .$$ Then the
set ${\ell}_{i,j}$ and the multiplicity with which each given
non-zero number appears among linking numbers, can be recovered
from $w_1(q)$. That is, $w_1(q)$ determines the polynomial
$\prod_{(i,j)}(x-{\ell}_{i,j})$ where the product is taken over
pairs $(i,j)$ with non-zero ${\ell}_{i,j}$.
\end{theorem}
\begin{proof}
The formula 3.2 (where $t^i=T_i$), or straightforward computation
using homotopy skein relations, gives us:
\begin{formulla}\label{4.4}
$w_1(q) = \sum_{(i,j)}q^{\sum_{(k,l)\neq (i,j)}2{\ell}_{k,l}}
\frac{q^{2{\ell}_{i,j}}-1}{q-q^{-1}} = -q^{2lk(L)}\sum_{(i,j)}
\frac{q^{-2{\ell}_{i,j}}-1}{q-q^{-1}}$
\end{formulla}
where $lk(L)= \sum_{(i,j)}{\ell}_{i,j}$. It is now an easy
exercise to see that
$\sum_{(i,j)}\frac{q^{-2{\ell}_{i,j}}-1}{q-q^{-1}}$ determines the
(unordered) sequence $\{{\ell}_{i,j}\}$. Notice that
$w_1(1)=lk(L)$ and $w_0(q)=q^{2lk(L)}$.
\end{proof}
Let us take for a moment a slightly more general point of view.
Assume that $a_1,a_2,...,a_k$ is a sequence of integers. Form from
the sequence the ``Young diagram" (positive numbers ``build" the part of
the diagram in the first quadrant and negative in the third), so
that rows of the diagram corresponds to numbers $\{a_i\}$, compare
Fig.4.5.

\ \\
\centerline{\psfig{figure=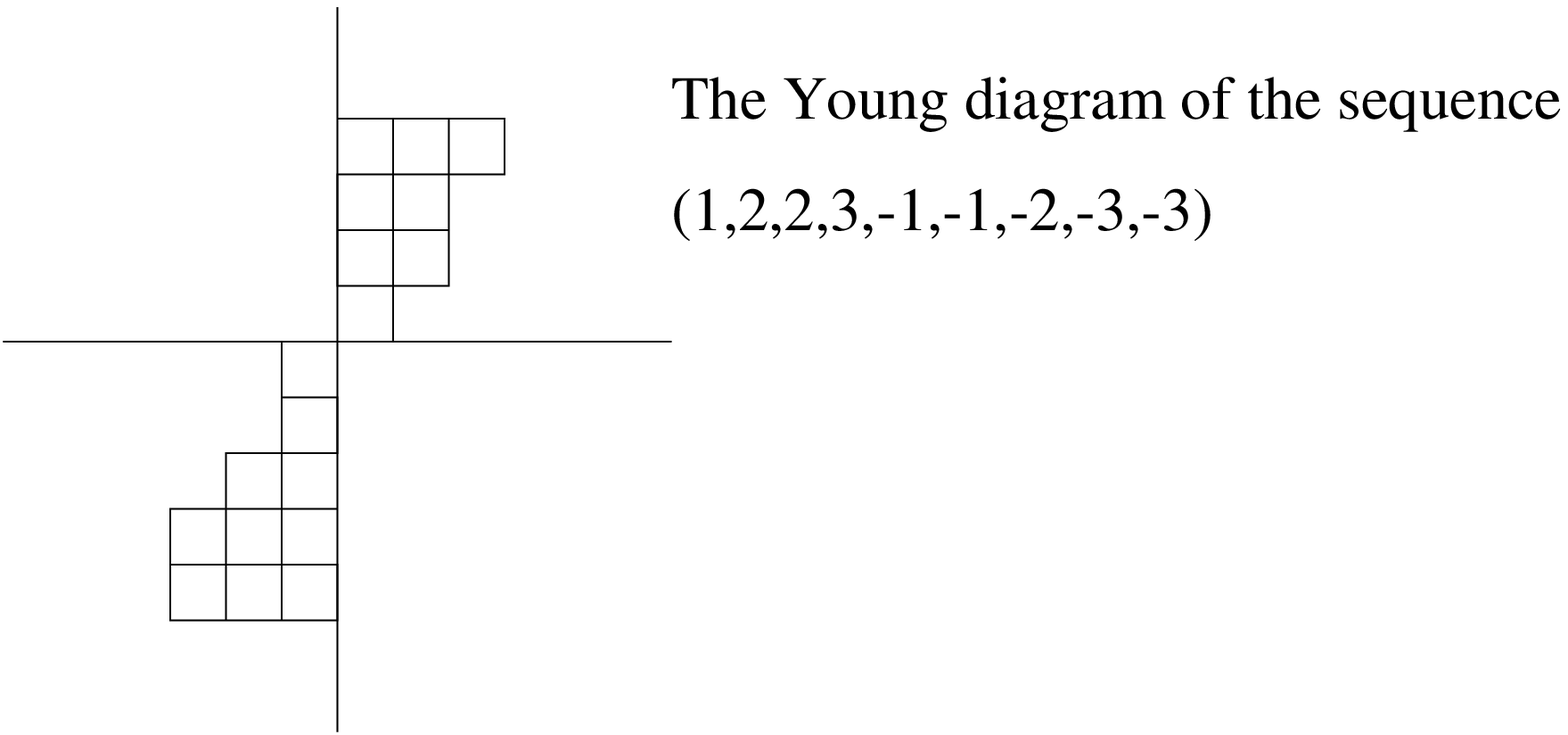,height=4.9cm}}
\begin{center}
Fig. 4.5.
\end{center}

\begin{lemma}\label{4.5}
The columns of the ``Young diagram" of the sequence
$a_1,a_2,...,a_k$ correspond to the coefficients of $\Sigma (q) =
\sum_{i=1}^k \frac{q^{2a_i}-1}{q-q^{-1}}$, that is, if \ $\Sigma
(q) = \sum_{j\neq 0} b_jq^{2j-sgn(j)}$ then $b_j$ is the number of
elements of the $j$'th column of the ``Young diagram" of the
sequence $\{a_i\}$, with appropriate signs. We can say shortly
that sequences $\{a_i\}$ and $\{b_j\}$ are dual one to another
(i.e. they represent dual ``Young diagrams").
\end{lemma}

\begin{proof}
It follows immediately from the identities:
$$\frac{q^{2a}-1}{q-q^{-1}}= q^{2a-1}+q^{2a-3}+...+q \ \ for \ \
a>0 ,$$ and $$\frac{q^{2a}-1}{q-q^{-1}}=
-q^{2a+1}-q^{2a+3}-...-q^{-1}\ \ for \ \ a<0,$$ which describe a
positive (resp. negative) row of length $a$ of the ``Young
diagram".
\end{proof}

\begin{corollary}\label{4.6}
The polynomial $w_1(q)$ has the form: $$w_1(q) =
-q^{2lk(L)}\sum_{j\neq 0} b_jq^{2j-sgn(j)}$$ where $b_j$'s are
heights of the columns of the ``Young diagram" corresponding to the
sequence $\{-\ell_{i,j}\}$
\end{corollary}

\begin{corollary}\label{4.7}
The polynomial $w_1(q)$ is unimodal, that is, if $w_L(q)= \sum
c_iq^{2i+1}$ then there is $j$ such that $...\leq |c_{j-1}| \leq
|c_{j}| \geq |c_{j+1}| \geq ...$.
\end{corollary}

\begin{example}\label{4.8}
For the link $L_2$ of Fig. 4.1 one has linking numbers $(2,3,-2)$.
Corollary 4.6 allows us to find immediately $w_1(q)$ by building
the Young diagram for numbers $(-2,-3,2)$ as in Fig. 4.6. \ \\
\centerline{\psfig{figure=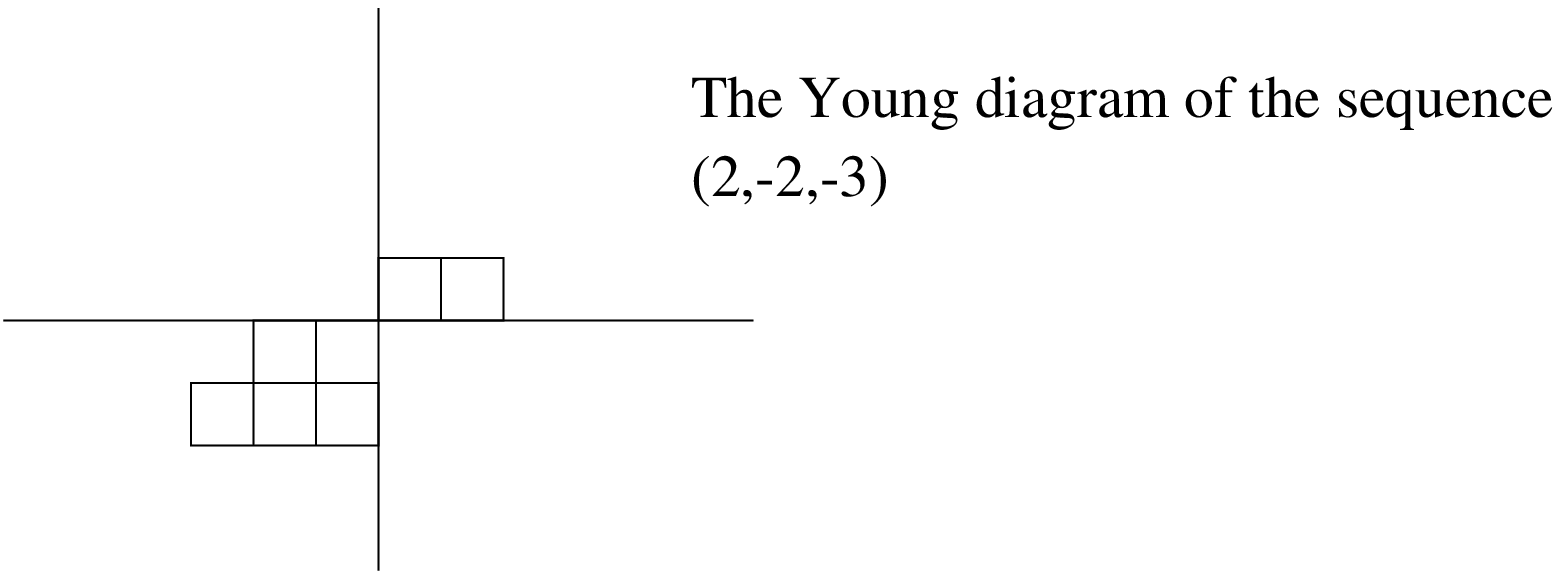,height=4.5cm}}
\begin{center}
Fig. 4.6.
\end{center}
\ \\ Thus, according to Corollary 4.6 one gets:
$$w_1(q)=-q^6(-2q^{-1} - 2q^{-3} - q^{-5} + q + q^3)= q+ 2q^3 +
2q^5 -q^7 - q^9$$ as in Example 4.1.
\end{example}

\begin{remark}\label{4.9}
It follows from Corollary 4.6 that non-zero linking numbers,
${\ell}_{i,j}$ can be recovered from $w_1(q)$ and, vice versa,
$w_1(q)$ can be recovered from the non-zero linking numbers.
Furthermore, the coefficient of $z$ in $HP_L(q,t,z)$ is equal to
$w_1(q)t^{n-1}$, so from this the number of components of $L$ can
be recovered (it is $n$), as well as the number of ${\ell}_{i,j}$
which are equal to $0$.
\end{remark}

\section{Elementary properties of homotopy skein modules}
It is convenient to consider homotopy and $q$-homotopy skein
modules, as special cases of homotopy skein modules with general
coefficients.
\begin{definition}\label{5.1}
We define ${\cal H S}(M;R,q,z)=R{\cal L}^h/S$ where $S$ is the
submodule generated by the expressions $q^{-1}L_{+}- qL_{-} -zL_{0}$
arising from mixed crossings, and $R$ is a commutative ring with identity,
$q$ an invertible element and $z$ any element of $R$.

One has of course ${\cal H S}(M;Z[z],1,z)={\cal H S}(M)$ and
${\cal H S}(M;Z[q^{\pm 1},z],q,z)= {\cal H S}^q(M)$.
\end{definition}
We list below a few useful elementary properties of homotopy skein
modules. Proofs of these properties are analogous to that of other
skein modules, \cite{P-S,Pr-7,Pr-8} and we omit them.
\begin{theorem}\label{5.2}
\begin{enumerate}
\item [(1)] An orientation preserving embedding of 3-manifolds
$i: M \to N$ yields a homomorphism of skein modules $i_{*}:
{\cal HS}(M;R,q,z) \to {\cal HS}(N;R,q,z)$. The above
correspondence leads to a functor from the category of 3-manifolds
and orientation preserving embeddings (up to ambient isotopy) to
the category of $R$-modules (with a specified element $z\in R$,
and an invertible element $q \in R$).
\item [(2)]
(Universal Coefficient Property)\\ Let $r: R \to R'$ be a
homomorphism of rings (commutative with 1). We can think of $R'$
as an $R$ module. Then the identity map on ${\cal L}$ induces the
isomorphism of $R'$ (and $R$) modules: $$ \bar r: {\cal
HS}(M;R,q,z)\otimes_{R} R' \to {\cal HS}(M;R',r(q),r(z)) .$$ In
particular ${\cal HS}(M;R,q,z)={\cal HS}^q \otimes_{Z[q^{\pm
1},z]}R$.
\item [(3)] Let $M=F\times I$ where $F$ is an oriented surface.
Then ${\cal H S}(M;R,q,z)$ is an algebra, where $L_{1} \cdot
L_{2}$ is obtained by placing $L_{1}$ above $L_{2}$ with respect
to the product structure. The empty link $T_{0}$ is the neutral
element of the multiplication. Every embedding $i:F' \to F$ yields
an algebra homomorphism $i_{*}: {\cal H S}(F'\times I;R,q,z) \to
{\cal H S}(F\times I;R,q,z)$.
\end{enumerate}
\end{theorem}
\section{The case of $M=F\times I$}\label{6}
Let ${\cal L}^{h}$ denote the set of homotopy links in $M$, that is,
${\cal L}^{h} = {\cal L}/(L_{+} - L_{-})$ where relations are
yielded by self-crossings. Let $\hat \pi$ denote
the set of conjugacy classes in $\pi_{1} (M)$, or equivalently the
set of homotopy knots in $M$. Choose some linear ordering, denoted
by $\leq$, of elements of $\hat \pi$. Given a homotopy link
$L=\{K_{1},K_{2},\ldots ,K_{n}\}$ in $F\times I$, we shall say
that $L$ is a layered homotopy link with respect to the ordering
of $\hat \pi$ if each $K_{i}$ is above $K_{i+1}$ in $F\times I$
and $K_{i}\leq K_{i+1}$. Let ${\cal B}$ be the set of all layered
homotopy links with respect to the ordering of $\hat \pi$,
including the empty link.
\begin{theorem}\label{6.1}
\begin{enumerate}
\item[(i)] The $q$-homotopy skein module ${\cal H \cal S}^{q}(F\times I)$
is generated by ${\cal B}$.
\item[(ii)] The homotopy skein module ${\cal H \cal S}(F\times I)$
is freely generated by ${\cal B}$; \cite{H-P-1}.
\item[(iii)] If $\pi_{1} (F)$ is abelian then the $q$-homotopy skein
module ${\cal H \cal S}^{q}(F\times I)$ is freely generated by
${\cal B}$.
\end{enumerate}
\end{theorem}
\begin{proof}
We will use the following notation: If $D$ is a link diagram in
$F\times I$ and $p_{1},\ldots ,p_{s}$ are some of its crossings
then $D^{p_{1},\ldots ,p_{s}}_{\epsilon_{1},\ldots ,\epsilon_{s}}$
denotes the link diagram obtained from $D$ by choosing at $p_{i}$
positive or negative crossing, or smoothing depending on whether
$\epsilon_{i}$ is equal to $+$ , $-$ or $0$. Consider a map
$\alpha : R{\cal B} $ to ${\cal H \cal S}^q(F\times I)$ given by
$\alpha (L) =L$. We will prove that $\alpha$ is an isomorphism if
either $q=\pm 1$ or $\pi_{1}(F)$ is abelian.

We follow closely the proof of Theorem 2.1 in \cite{H-P-1} and the
proof of Theorem 1.6 in \cite{P-T}.

We will construct the inverse map $W$,
 to $\alpha$. The plan for
constructing $W$ is as follows:
\begin{enumerate}
\item[(i)]
We define $W$ on the diagrams by inducting on the number of
components and number of ``bad" crossings.
\item[(ii)]
The initial definition, in the inductive step, depends on the
ordering of components of a diagram and on the order of
eliminating bad crossings. We prove that our choices do not
give different results if either $q=\pm 1$ or $\pi_1(F)$ is
abelian.
\item[(iii)]
We show that $W$ is invariant under Reidemeister moves, and satisfies
the homotopy skein relations.
\end{enumerate}

To define $W$ we first we use induction on
the number of components, $c(D)$, of the link diagram $D$. For each
$n\geq 0$ we define a function $W_n$ defined on the set of
oriented link diagrams with no more than $n$ components. Then $W$
will be defined for every diagram by $W(D)=W_n(D)$ where $n\geq
c(D)$. Of course the functions $W_n$ must satisfy certain
coherence conditions for this to work. First we put
$W_0(\emptyset) = \emptyset$ and $W_1(D_K)= K$ where $D_K$ is a
diagram of a knot $K$.

To define $W_{n+1}$
and prove its properties we will use induction several times. The
following will be called the ``Main Inductive Hypothesis":
M.I.H. We assume that we have already defined a function $W_n$ for
each diagram $D$ with no more than $n$ components ($c(D)\leq n$).
We assume that $W_n$ has the following properties:
\begin{enumerate}
\item[(1)] $W_n(D)=L_D$ if $D$ is a layered diagram of $n$ or less components
representing a layered link ($L_D$) respecting the $\leq$ ordering of
$\hat\pi$.
\item[(2)] $W_n(D_+) = W_n(D_-)$ for a self-crossing ($c(D_+)\leq n$),
\item[(3)]
$q^{-1}W_n(D_+) - qW_n(D_-)=zW_n(D_0)$ for a mixed crossing
($c(D_+)\leq n$).
\item[(4)] $W_n(D^R)=W_n(D)$ where $D^R$ is the result of
a Reidemeister move on $D$ ($c(D)\leq n$).
\end{enumerate}
Then we want to make the Main Inductive Step, M.I.S., to obtain
the existence of a function $W_{n+1}$ with analogous properties
defined on diagrams with at most $n+1$ components.

Before dealing with the task of making the M.I.S. let us explain
that it will end the proof of the theorem. It is clear that the
function $W_{n}$ satisfying M.I.H. is uniquely determined by
properties (1)-(3) and the fact that any diagram can be changed to
a layered diagram (respecting the ordering $\leq$ on $\hat \pi$) 
by changing some
crossings (and observing that smoothing is lowering the number of
components). Thus the compatibility of the functions $W_n$ is
obvious and they define a function $W$ on diagrams. The function
satisfies skein relations by (2)-(3) and Reidemeister moves by
(4). By property (1) it is the inverse function to $\alpha$.

The rest of the section will be occupied by M.I.S.

First we define a function $W_d$ on diagrams having $n+1$
components ($n\geq 1$) where components are ordered according to their
homotopy type and the chosen ordering of $\hat\pi$. If no two of
the components are of the same homotopy type then their ordering
is uniquely determined. However, if some of the components have
the same homotopy type then any choice when ordering them is
allowed. We proceed by induction on the number, $b(D)$, of bad
crossings of an ordered diagram $D$, that is, crossings that have to
be changed in order to obtain a diagram of a layered link (in the
chosen ordering of components). If $b(D)=0$ then we put
$W_d(D)=L_D$ where $L_D$ is a layered link represented by $D$.

Assume that $W_d$ is defined for all $D$ with $b(D)<t$, $t\geq 1$.
Let $D$ be an ordered $n+1$-component link diagram with $b(D)=t$.
In the next steps, (a)-(c), we define $W_d(D)$ and prove its
properties.
\begin{enumerate}
\item[(a)]\ \\
Let $p$ be a bad crossing of $D$. We define a function $W_{(p)}$
by the formula (depending on $p$ being positive or negative we
have $D=D^p_+$ or $D=D^p_-$):

$$ W_{(p)}(D)= \left \{
\begin{array}{ll} q^2W_d(D^p_-) + qz W_n(D^p_0) & \mbox{if
sgn$p=+$} \\ q^{-2}W_d(D^p_+) - q^{-1}z W_n(D^p_0) & \mbox{if
sgn$p=-$}
\end{array}
\right. $$

The right-hand side of the equation is defined by inductive
assumptions. We show that $W_d$ is independent on the choice of
a bad crossing of $D$. Let $s$ be another bad crossing and $
W_{(s)}(D)$ associated with $s$ value. We show that 
$W_{(s)}(D)=W_{(p)}(D)$. 

If $p$ and $s$ are mixed crossings between different pairs of 
components then we get $W_{(s)}(D)=W_{(p)}(D)$ without restricting
$q$ or the fundamental group. The computation is the same as for the 
Jones-Conway polynomial (see \cite{P-T}), 
so we perform this only in the case of $sgn(p) = sgn(s) = 1$:
$$W_{(p)}(D)=q^2W_d(D^p_-) + qz W_n(D^p_0) =$$
$$q^4W_d(D_{-,-}^{p,s}) + q^3zW_n(D_{-,0}^{p,s}) +
q^3zW_n(D_{0,-}^{p,s}) + q^2z^2W_n(D_{0,0}^{p,s}).$$
The result is symmetric with respect to $p$ and $s$, thus
$W_{(s)}(D)=W_{(p)}(D)$.

Now assume that $p$ and $s$ are mixed crossings between the same
pair of components of $D$.

For $q=\pm 1$ we will check all sign cases at once. 
Let $\epsilon (x)$ be the sign of a crossing $x$, and
$\epsilon (q)$ the sign of $q$. Then we have:
$$W_{(p)}(D)= W_d(D^{\ \ \ p}_{-\epsilon (p)}) + \epsilon
(p)\epsilon (q)zW_n(D^p_0)=$$ $$W_d(D^{\ \ \ p,\ \ \ s}_{-\epsilon
(p),-\epsilon (s)}) + \epsilon (s)\epsilon (q)zW_n(D^{\ \ p,\ \
s}_{-\epsilon (p),0}) +\epsilon (p)\epsilon (q)zW_n(D^p_0),$$
$$W_{(s)}(D)= W_d(D^{\ \ \ s}_{-\epsilon (s)}) + \epsilon
(s)\epsilon (q)zW_n(D^s_0)=$$ $$W_d(D^{\ \ \ s,\ \ \ p}_{-\epsilon
(s),-\epsilon (p)}) + \epsilon (p)\epsilon (q)zW_n(D^{\ \ \ s, \
p}_{-\epsilon (s),0}) +\epsilon (s)\epsilon (q)zW_n(D^s_0).$$

 By
(2) of the M.I.H. we have $W_n(D^{\ \ \ p,\ s}_{-\epsilon
(p),0})= W_n(D^s_0)$, and $W_n(D^p_0)=W_n(D^{\ \ \ s, \
p}_{-\epsilon (s),0})$. Thus $ W_{(s)}(D)=W_{(p)}(D)$ and we put
$W_d(D)=W_{(p)}(D)$, independent of the choice of a bad
crossing.

In the case when $\pi_1{(F)}$ is abelian but there are no
restrictions on $q$ we will check positive, negative and mixed
cases of signs of $p$ and $s$ separately.
\begin{enumerate}
\item[(++)] If $sgn(p)=sgn(s)= +1$ we get:
$$W_{(p)}(D)= q^2W_d(D^{p}_{-}) + qz
W_n(D^p_0)=q^4W_d(D^{p,\ s}_{-,-}) + q^3zW_n(D^{p,s}_{-,0})
+qzW_n(D^p_0),$$ 
$$W_{(s)}(D)= q^2W_d(D^{s}_{-}) + qz
W_n(D^s_0)=q^4W_d(D^{s,\ p}_{-,-}) + q^3zW_n(D^{s,p}_{-,0})
+qzW_n(D^s_0),$$
\item[(-- --)] If $sgn(p)=sgn(s)=-1$ we get:
 $$W_{(p)}(D)= q^{-2}W_d(D^{p}_{+}) - q^{-1}zW_n(D^p_0)=$$
$$ q^{-4}W_d(D^{p,\ s}_{+,+}) - q^{-3}zW_n(D^{p,s}_{+,0})
-q^{-1}zW_n(D^p_0),$$ 
$$W_{(s)}(D)= 
q^{-2}W_d(D^{s}_{+}) - q^{-1}z W_n(D^s_0)=$$
$$ q^{-4}W_d(D^{s,\ p}_{+,+}) - q^{-3}zW_n(D^{s,p}_{+,0})
-q^{-1}zW_n(D^s_0),$$
\item[(+--)] If $sgn(p)=+1$ and $sgn(s)=-1$ we get:
$$W_{(p)}(D)= q^2W_d(D^{p}_{-}) + qz W_n(D^p_0)=
W_d(D^{p,\ s}_{-,+})- qz W_n(D^{p,s}_{-,0}) + qz W_n(D^p_0),$$
$$W_{(s)}(D)= q^{-2}W_d(D^{s}_{+}) - q^{-1}z W_n(D^s_0)=$$
$$W_d(D^{s,\ p}_{+,-})+ q^{-1}z W_n(D^{s,p}_{+,0}) - q^{-1}z W_n(D^s_0),$$
\item[(--+)] The same as (+--) with the role of $p$ and $s$ switched.
\end{enumerate}
By (2) of the M.I.H. we have $W_n(D^{\ \ \ p,s}_{-\epsilon
(p),0})= W_n(D^s_0)$, and $W_n(D^p_0)=W_n(D^{\ \ \ s,p}_{-\epsilon
(s),0})$. Furthermore, for $\pi_1(F)$ abelian we have $W_n(D_0^p)=
W_n(D_0^s)$ (for $n=2$ it is immediate as $D_0^p$ and $D_0^s$ are 
homotopic knots; in general one can use induction on the number
of crossing at which $D_0^p$ is below the rest of the diagram).
Thus $ W_{(s)}(D)=W_{(p)}(D)$ and we put
$W_d(D)=W_{(p)}(D)$, independent of the choice of a bad
crossing. We should stress here that in cases (+--) and (--+) 
we have got the equality $W_d(D^{\ p,\ \ s}_{\epsilon(p),\epsilon(s)}) =
W_d(D^{\ \ p,\ \ s}_{-\epsilon(p),-\epsilon(s)})$.

\item[(b)] Homotopy skein relations.\\
The fact that $W_{d}$ satisfies homotopy skein relation follows
from the construction for the mixed crossing and by an easy
induction on $b(D)$ for a self-crossing:
\begin{enumerate}
\item[(m.c.)] Let $p$ be a mixed crossing of $D$. Then $p$ is a bad
crossing of $D^p_+$ or $D^p_-$. Using the defining relation for
$W_{(p)}$ for the diagram in which $p$ is a bad crossing we get
the required skein relation: $q^{-1}W_d(D^p_+) - qW_d(D^p_-)
=zW_n(D^p_0)$.
\item[(s.c.)] Let $p$ be a self-crossing of $D$. We proceed
by induction on $b(D)$ of the number of bad crossings of an
$n+1$ component ordered link diagram $D$. If $b(D)=0$ then
$b(D^p_{-sgn(p)})=0$ as well and $W_{d}(D^p_+)=W_{d}(D^p_-) = L_D$
the layered link represented by $D$. Assume the skein relation
holds for $b(D)<t$, $(t>0)$. Assume $b(D)=t$ and let $s$ be its
bad crossing. Then $W_{d}(D^p_+)= W_{d}(D^{p,\ s}_{+,\epsilon
(s)})= q^{2 \epsilon (s)} W_{d}(D^{p,\ s}_{+,-\epsilon (s)}) +
\epsilon (s) q^{\epsilon (s)} W_{n}(D^{p,s}_{+,0}) = q^{2 \epsilon (s)}
 W_{d}(D^{p,\ s}_{-,-\epsilon (s)}) +$\\ 
$\epsilon (s) q^{\epsilon (s)}W_{n}(D^{p,s}_{-,0}) = W_{d}(D^p_-)$, 
by inductive assumptions.
\end{enumerate}
\item[(c)] Independence of $W_d$ on Reidemeister moves.
\begin{enumerate}
\item[($\Omega 1$)]
Let $\Omega_1 (D)$ be obtained from $D$ by the first Reidemeister
move, where $D$ is an ordered link diagram of $n+1$ components. We
proceed by induction on $b(D)$. We have $b(\Omega_1(D))=b(D)$ and
if $b(D)=0$ then $W_{d}(D)=L_D =W_{d}(\Omega_1(D))$. Assume that
equality holds for $b(D)<t$, ($t>0$), and let $b(D)=t$. Finally
let $p$ be a bad crossing of $D$ (and $\Omega_1(D)$). Then
$W_{d}(D)= q^{2\epsilon (p)}W_{d}(D^p_{-\epsilon (p)}) + \epsilon
(p)q^{\epsilon (p)} W_{n}(D^p_{0})= q^{2\epsilon
(p)}W_{d}(\Omega_1(D)^p_{-\epsilon (p)}) + \epsilon (p)q^{\epsilon
(p)}W_{n}(\Omega_1(D)^p_0) = W_{d}(\Omega_1(D))$, by inductive
assumptions.
\item[($\Omega 2$)]
Let $\Omega_2(D)$ be obtained from $D$ by a second Reidemeister
move, where $D$ is an ordered link diagram of $n+1$ components. We
have to consider two cases: Reidemeister move does not create any
new bad crossings or it creates two new bad crossings.
\begin{enumerate}
\item[(i)] In the first case we proceed by induction on $b(D)$,
exactly as in the case of $\Omega_1(D)$.
\item[(ii)] $\Omega_2(D)$ is introducing two new bad crossings of
opposing signs $p$ and $s$ (we can assume $sgn(p)=+$ and $sgn
(s)=-$) (see Fig. 6.1). In particular $p$ and $s$ are mixed
crossings. Then:\\
\ \\
\centerline{\psfig{figure=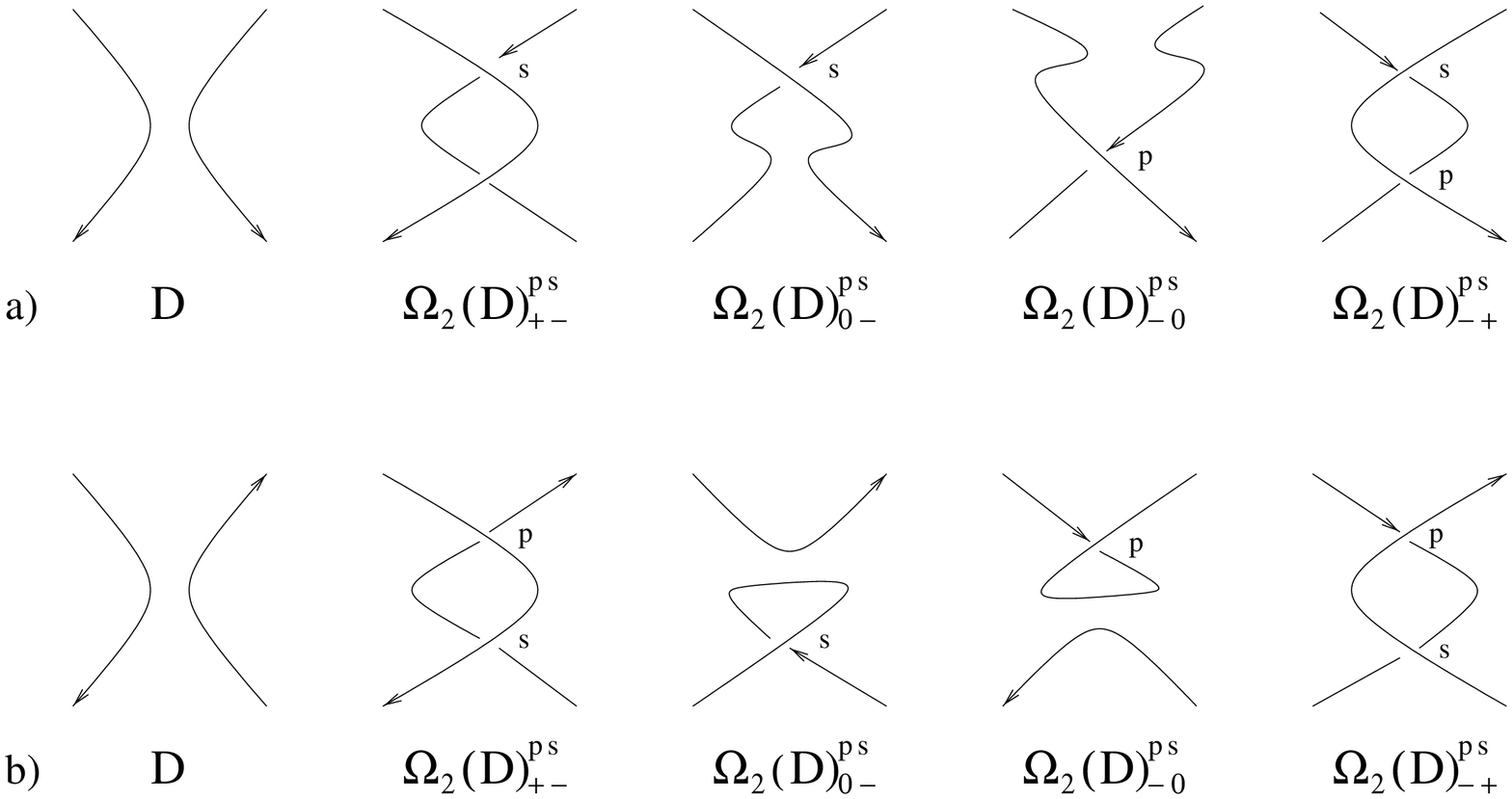,height=6.2cm}}
\begin{center}
Figure 6.1
\end{center}
\ \\
$W_{d}(\Omega_2(D)) = W_{d}(\Omega_2(D)^{p,s}_{+,-})=
W_{(p)}(\Omega_2(D)^{p,s}_{+,-})=q^{2}W_{d}(\Omega_2(D)^{p,s}_{-,-})
+ qzW_{n}(\Omega_2(D)^{p,s}_{0,-})= W_{d}(\Omega_2(D)^{p,s}_{-,+})
+qz (W_{n}(\Omega_2(D)^{p,s}_{0,-}) -
W_{n}(\Omega_2(D)^{p,s}_{-,0})) = W_{d}(\Omega_2(D)^{p,s}_{-,+})=
W_{d}(D)$, by inductive assumptions. More precisely:
$W_{d}(\Omega_2(D)^{p,s}_{-,+})= W_{d}(D)$ as
$\Omega_2(D)^{p,s}_{-,+}$ is obtained from $D$ by a second
Reidemeister move preserving $b(D)$, and either
$\Omega_2(D)^{p,s}_{0,-} = \Omega_2(D)^{p,s}_{-,0}$ (Fig 6.1 (a))
or they differ by first Reidemeister moves (Fig 6.1 (b)).
\end{enumerate}
\item[($\Omega 3$)] Let $p_1$ be the top crossing of the third 
Reidemeister move,
$p_3$ be a bottom crossing and $p_2$ be a crossing between the top
and bottom arcs of the move (see Fig.6.2). Notice that
$b(\Omega_3(D))=b(D)$. To show that $W_d(\Omega_3(D))= W_b(D)$,
we proceed by induction on $b(D)$. If $b(D)=0$ then $W_{d}(D)=L_D=
W_{d}(\Omega_3(D))$. Assume that the equality
$W_{d}(D)=W_{d}(\Omega_3(D))$ holds for $b(D)<t$ ($t>0$), and let
$b(D)=t$. Then either:\\
 \ \\
\centerline{\psfig{figure=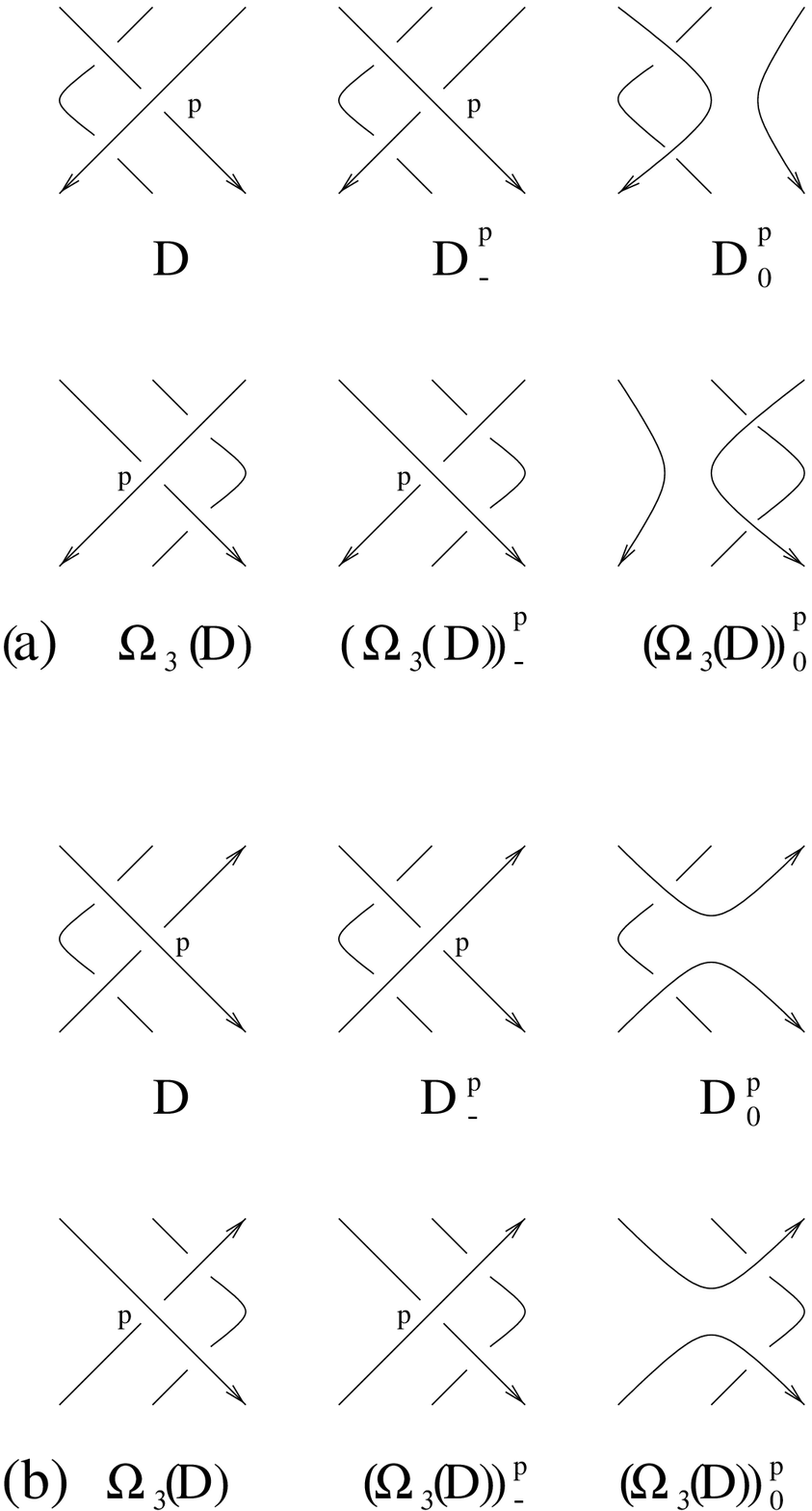,height=10.8cm}}
\begin{center}
Figure 6.2.
\end{center}

\begin{enumerate}
\item[(i)] there is a bad crossing $p$ different than $p_i$ in $D$.
Then: $W_{d}(D)= q^{2\epsilon (p)}W_{d}(D^{\ p}_{-\epsilon (p)}) +
\epsilon (p)q^{\epsilon (p)} W_{n}(D^p_{0}) = $\\
$q^{2\epsilon
(p)}W_{d}(\Omega_3(D)^p_{-\epsilon (p)}) + \epsilon (p)q^{\epsilon
(p)}W_{n}(\Omega_3(D)^p_0) = W_{d}(\Omega_3(D))$, by inductive
assumptions.
\item[(ii)] Assume that $p=p_1$ or $p_3$ is a bad crossing (assume
for simplicity of notation that $sgn(p)=+$). Then:
$W_{d}(\Omega_3(D)) = q^2W_{d}(\Omega_3(D)^p_-) + qz
W_{n}(\Omega_3(D)^p_0) = q^2W_{d}(\Omega_3(D^p_-)) + qzW_n(\Omega_3(D)^p_0)=
q^2W_d(D^p_-) + qzW_n(D^p_0) = W_d(D)$, by inductive assumptions
(compare Fig. 6.2). More precisely: 
because $b(D^p_-)< b(D)$
therefore by the inductive assumption 
$W_d(\Omega_3(D)^p_-) =W_d( \Omega_3(D^p_-))= W_d(D^p_-)$. 
Furthermore, by M.I.H. we have
$W_n(\Omega_3(D)^p_0) =W_n( \Omega(D^p_0))= W_n(D^p_0)$, where
$\Omega$ can be a composition of second Reidemeister moves (Fig.6.2(a))
or ambient isotopy (Fig.6.2(b)).
\item[(iii)] Assume that $p_2$ is a bad crossing then from the definition
of $p_1,p_2$ and $p_3$ and from the fact that we perform a third
Reidemeister move, follows that $p_1$ or $p_3$ is a bad
crossing, so we deal with the case (ii).
\end{enumerate}
\item[(c)] Independence of $W_d$ on the ordering of components.\\
Our first choice, when defining $W_d$, was the chosen ordering of
components of a diagram. As the choice was coherent with the fixed
ordering of $\hat\pi$, the only freedom was in the choice of
ordering of components having the same homotopy type. Consider two
orderings, $D_1$ and $D_2$, of the given $n+1$ component 
layered link diagram
$D$, which differ only by switching the order of two neighboring
components $D'$ and $D''$ which are homotopic (in $D_1$
one has $D''$ before $D'$ and in $D_2$ one has $D'$ before $D''$).
To prove independence of $W_d$ on an ordering of homotopically equal
components it suffices to show that $W_d(D_1)=W_d(D_2)$.

First change $D'$ in $D_1$ by Reidemeister moves so that $D'$ and
$D''$ are parallel with $D''$ on top slightly to the ``left" of $D'$
(compare Fig. 6.3). By (b), the value of $W_d$ is unchanged. Now,
crossings between $D'$ and $D''$ are grouped in pairs of crossings
of opposite signs; $c^+$ and $c^-$ of Fig.6.3. Using homotopy
skein relations to crossings $c^+$ and $c^-$, we get:\\ 
$ W_d(D_1)
= W_d((D_1)^{c^+,c^-}_{+,-}) = q^2W_d((D_1)^{c^+,c^-}_{-,-}) +
qzW_n((D_1)^{c^+,c^-}_{0,\ \ -}) =
q^2(q^{-2}W_d((D_1)^{c^+,c^-}_{-,\ +}) -
q^{-1}zW_n((D_1)^{c^+,c^-}_{-,\ 0})) + qzW_n((D_1)^{c^+,c^-}_{0,\ -})
= W_d((D_1)^{c^+,c^-}_{-,\ +}) + qz(W_n((D_1)^{c^+,c^-}_{0,\ -}) -
W_n((D_1)^{c^+,c^-}_{-,\ 0})) = W_d((D_1)^{c^+,c^-}_{-,+}) $. 
We use here the fact that $(D_1)^{c^+,c^-}_{0,\ -}$ and
$(D_1)^{c^+,c^-}_{-,\ 0}$ are homotopic layered links of $n$ components 
($(D''\cup D')^{c^+,c^-}_{0,\ -}$ and $(D''\cup D')^{c^+,c^-}_{-,\ 0}$ are
homotopic knots). Applying our identity to any pair of crossings
between $D''$ and $D'$ one gets $W_d(D_1)=W_d(D_2)$. Thus $W_d$
does not depend on an ordering of homotopic components.
\end{enumerate}
(a)-(c) allows us to put $W_{n+1}(D)=W_d(D)$ for a link diagram
with $n+1$ components, and complete the Main Inductive Step.
\end{enumerate}
By construction $W$ is the two-sided inverse of $\alpha$ so the
proof of Theorem 6.1 is complete.
\end{proof}
\ \\ \ \\
\centerline{\psfig{figure=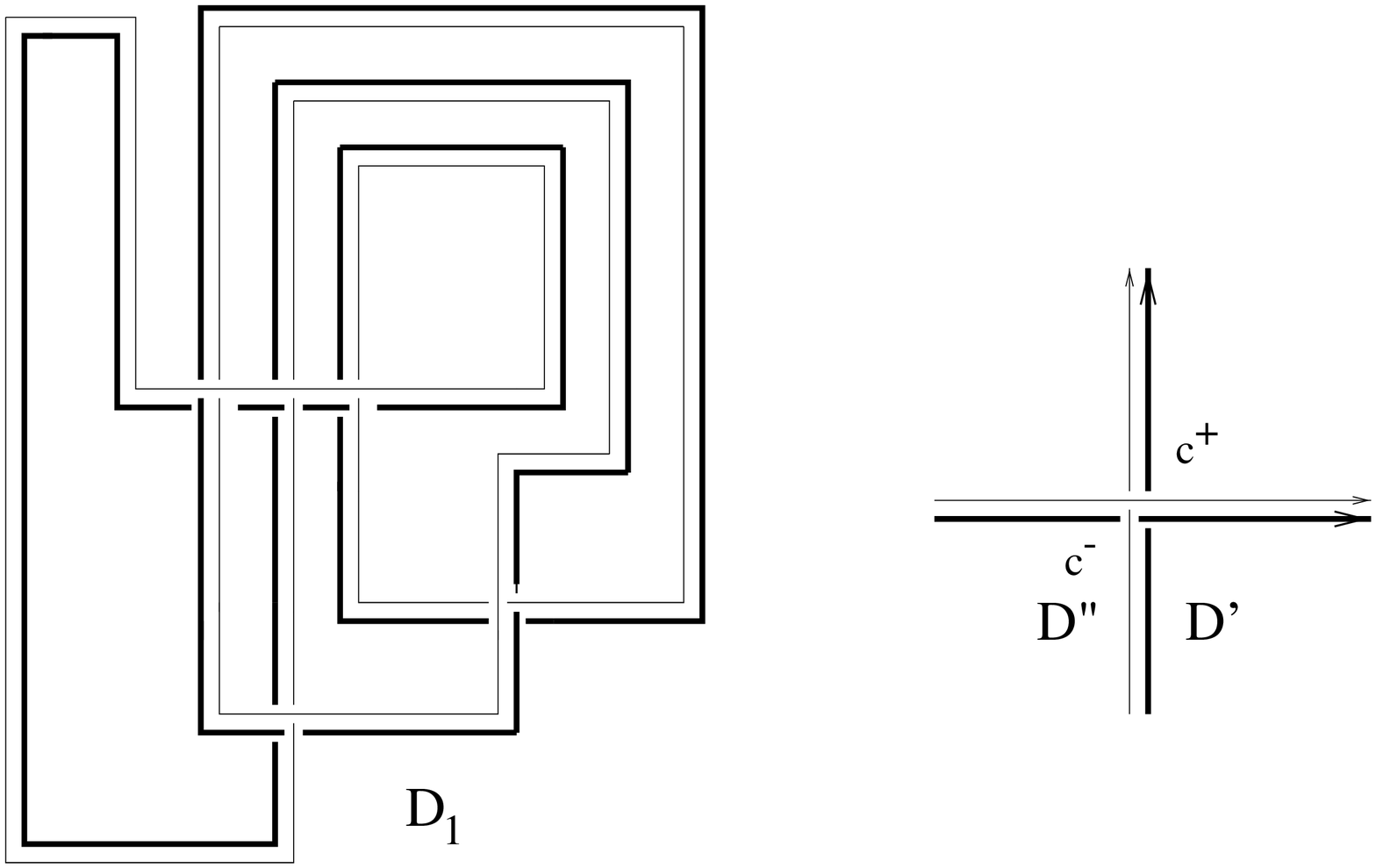,height=5.2cm}}
\begin{center}
Figure 6.3.
\end{center}
\ \\
\begin{remark}
Our assumption that $q= \pm 1$ or $\pi_1(F)$ is abelian was used
only in the step (a), where we proved independence of $W_d$ on the
ordering of ``bad" crossings. Careful analysis of this step would
allow us to identify the kernel of $\alpha : R{\cal B} \to {\cal H
\cal S}^q(F\times I)$, however I am not sure how to find good
description of it (one approach, reduction to homologies, is
sketched in Section 8). We will show, in the next section, that
$ker(\alpha)$ is not trivial for any $F$ with negative Euler
characteristic.
\end{remark}

\section{Torsion}\label{7}
If the Euler characteristic, $\chi (F)$, is negative then the
$q$-homotopy skein module has torsion and $ker(\alpha)\neq
\{0\}$. This is described in Theorem 7.1. On the other hand we
have proven in Theorem 6.1(ii) that for $q=\pm 1$ the module is free.
We have another very interesting case of coefficients, ${\cal
F}(q)[z]$, where ${\cal F}(q)$ is the field of rational functions
in variable $q$. We discuss this in Problem 7.3.

\begin{theorem}\label{7.1}
Let $F$ be a surface (not necessary compact) which contains a disc
with 2 holes or a torus with a hole embedded $\pi_1$-injectively;
equivalently, $\pi_1 (F_{0})$ is not abelian for a connected
component $F_{0}$ of $F$ (in the compact connected case this 
 means that $\chi (F)<0$). Then
\begin{enumerate}
\item[(a)]
${\cal H \cal S}^q(F\times I)$ has torsion.
\item[(b)] Let $\alpha : R{\cal B} \to {\cal H \cal S}^q(F\times I)$
be an $R$-homomorphism introduced in the previous section. Then
$ker \alpha \neq \{0\}$.
\end{enumerate}
\end{theorem}
\begin{proof}
Let $i: F' \to F$ be an embedding of surfaces which is
$\pi_1$-injective (that is, $i_*: \pi_1(F') \to \pi_1(F)$ is a
monomorphism), then by Theorem 6.1 (ii), the embedding yields a
monomorphism of homotopy skein modules, $i_{\#}:{\cal H
S}(F'\times I) \to {\cal H S}(F\times I)$. Now the theorem follows
from the following lemma, which describes torsion elements in a disc
with two holes and a torus with a hole, as torsion elements go to
torsion elements under the homomorphism $i_q:{\cal H
S}^q(F'\times I) \to {\cal H S}^q(F\times I)$.
Furthermore, torsion elements constructed in Lemma 7.2 are also
nontrivial for $q=1$ so they ``survive" the homomorphism $i_q$
by Theorem 6.1 (ii) and Theorem 5.2(2).
\end{proof}
\begin{lemma}\label{7.2}
\begin{enumerate}
\item[(1)] Let $F$ be a disk with two holes with $\pi_1(F)=\{x,y | \ \}$.
Consider the word $w=xy^{-1}$ and its inverse $\bar w = yx^{-1}$.
Let $\gamma$ and $\bar\gamma$ be two knots realizing these words
and $L=\gamma \bullet \bar{\gamma}$ be a two component link in
$F\times I$, being a product of these knots and $D_L$ the diagram
of $L$ with two mixed crossings, $p$ and $s$, as shown in Fig.
7.1. Then we can ``compute" $D_L$ in ${\cal H \cal S}^q(F\times
I)$ resolving first $p$ or $s$:
\begin{enumerate}
\item[(p,s)]
$D_L = q^2(D_L)^p_- + qz(D_L)^p_0 = (D_L)^{p,s}_{-,+} -
qz(D_L)^{p,s}_{-,0} + qz(D_L)^p_0 = (D_L)^{p,s}_{-,+} +
qz((D_L)^p_0 - (D_L)^{p,s}_{-,0})$.
\item[(s,p)]
$D_L = q^{-2}(D_L)^s_+ - q^{-1}z(D_L)^s_0 = (D_L)^{s,p}_{+,-} +
q^{-1}z(D_L)^{s,p}_{+,0}- q^{-1}z(D_L)^s_0= (D_L)^{p,s}_{-,+} +
q^{-1}z((D_L)^{s,p}_{+,0} - (D_L)^s_0)$.
\end{enumerate}
Thus in ${\cal H \cal S}^q(F\times I)$ one has:
$$(q-q^{-1})z((D_L)^p_0 - (D_L)^s_0) = 0$$ To see that we really
identified a nontrivial torsion element one should show that $(D_L)^p_0 -
(D_L)^s_0 \neq 0$, but this is the case even for $q=1$, as
$(D_L)^p_0$ is a knot representing $xy^{-1}x^{-1}y$ in $\pi_1(F)$
and $(D_L)^s_0$ is a knot representing $xyx^{-1}y^{-1}$. These
two elements are not conjugate in $\pi_1(F)$ so, by Theorem 6.1,
are different in the homotopy skein module of $F\times I$. \ \\ \
\\
\centerline{\psfig{figure=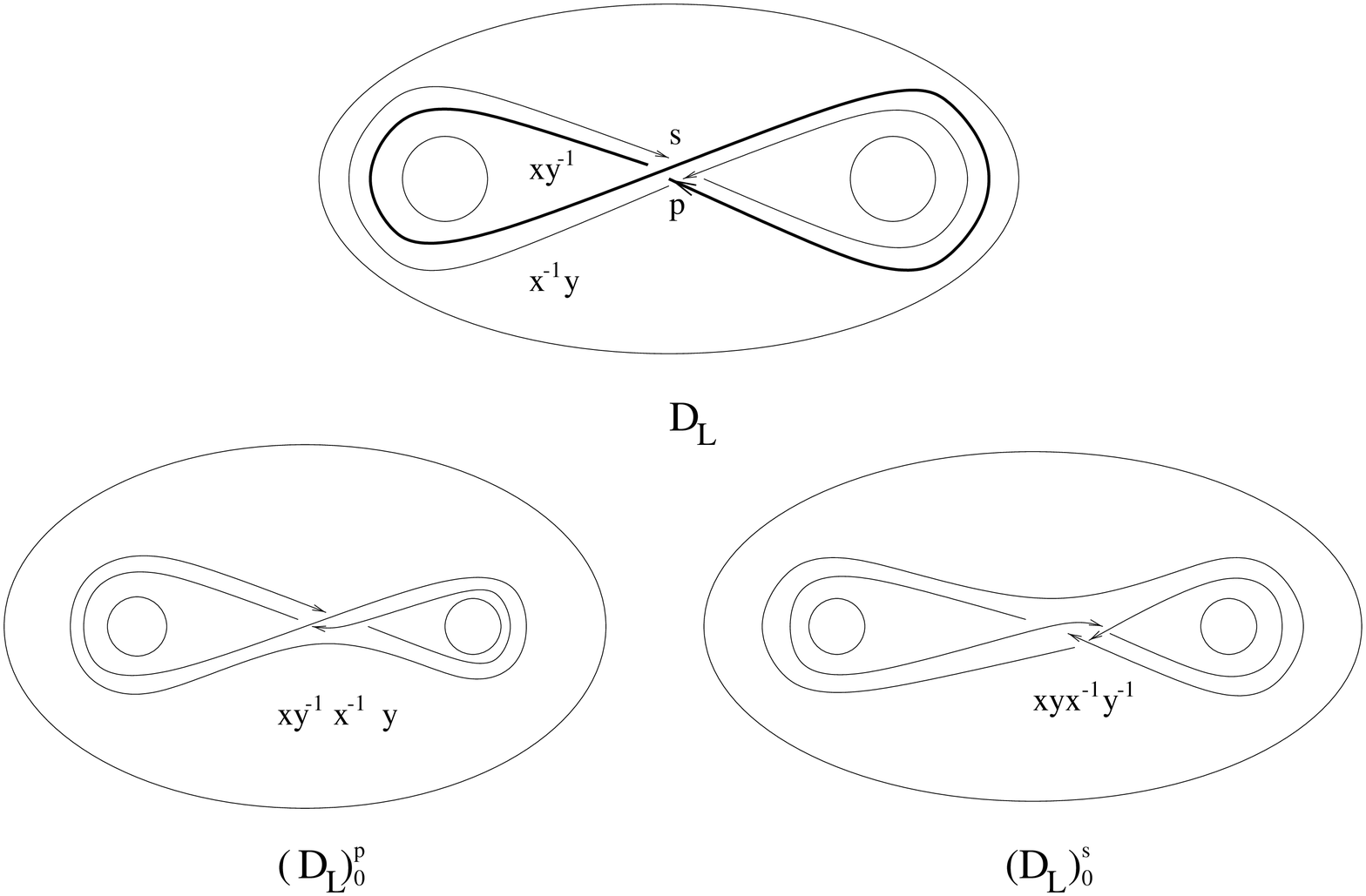,height=7.6cm}}
\begin{center}
Figure 7.1
\end{center}
\ \\

\item[(2)]
Let $F$ be a torus with a hole with $\pi_1(F)=\{x,y | \ \}$.
Consider the word $w=xy$ and the word $ w' = xy^{-1}$. Let $\gamma$
and ${\gamma}'$ be two knots realizing these words and $L=\gamma
\bullet {\gamma}'$ be a two component link in $F\times I$, being a
product of these knots. Let $D_L$ be the diagram of $L$ with two
mixed crossings, $p$ and $s$, as shown in Fig. 7.2. Then we can
``compute" $D_L$ in ${\cal H \cal S}^q(F\times I)$ resolving first
$p$ or $s$ and getting (similarly as in (1)): 
$(q-q^{-1})z((D_L)^p_0 - (D_L)^s_0) = 0$. 
Again we deal with a nontrivial torsion element because
elements of $\pi_1(F)$ represented by $(D_L)^p_0$ and $(D_L)^s_0$
are not conjugate.
\end{enumerate}
\end{lemma}
\ \\
\centerline{\psfig{figure=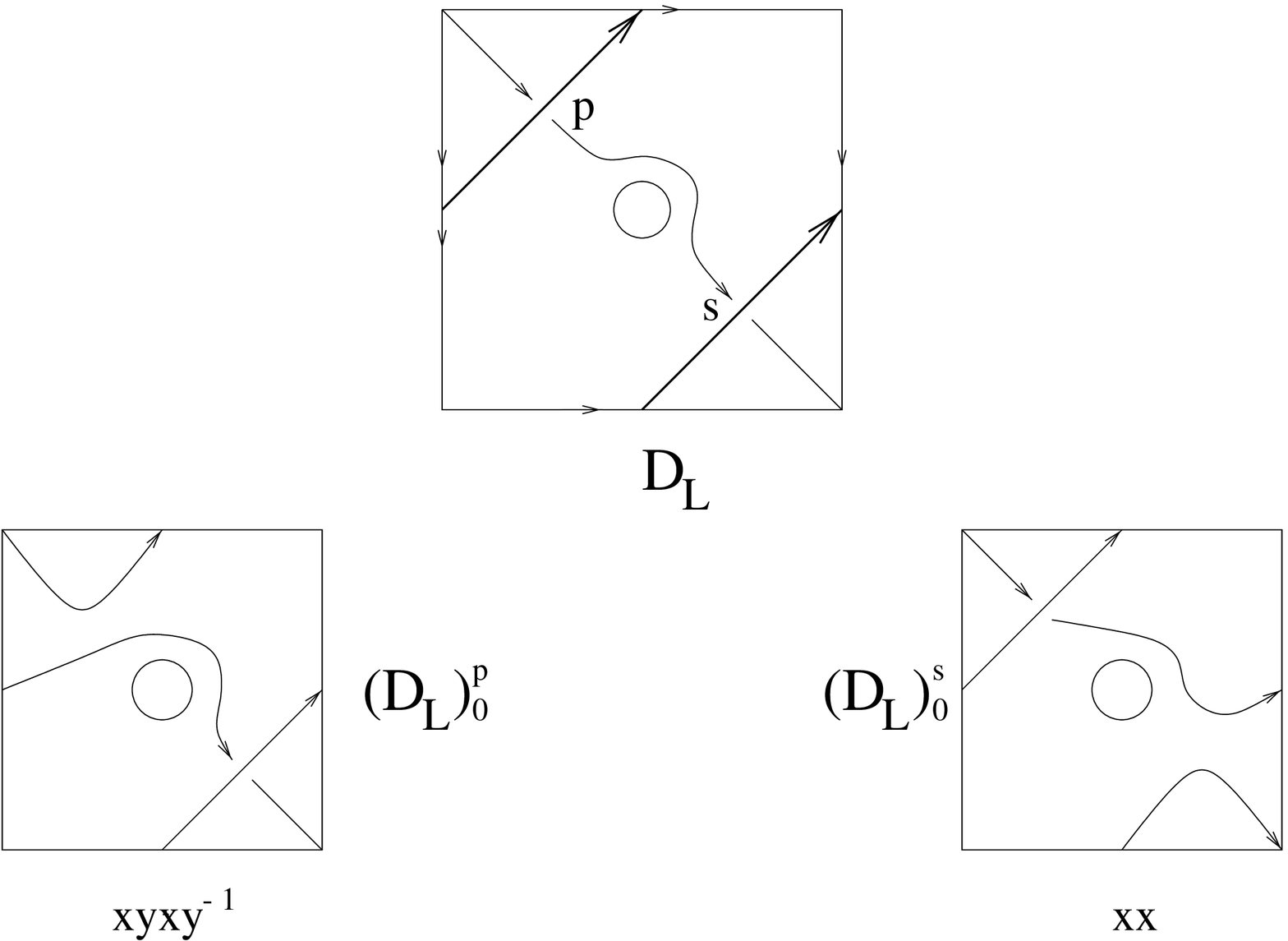,height=7.8cm}}
\begin{center}
Figure 7.2
\end{center}
\ \\

If we allow polynomials $q^{n}-1$ to be invertible in the ring of
coefficients, then our examples are not producing torsion but
instead reduce the number of generators in the skein module.

\begin{problem}
If two knots are homologous in $F\times I$, are they equal in the
homotopy skein module with coefficients in ${\cal F}(q)[z]$?

If the answer is yes then ${\cal HS}(F\times I;{\cal
F}(q)[z],q,z)$ is algebra isomorphic to the quantization of a
$q$-symmetric Poisson algebra described in Section 8 (see Theorem
8.11).
\end{problem}

\newpage

\section{Lie algebras, Poisson algebras, Universal enveloping algebras,
$q$-algebras, quantizations and relations to homotopy skein
modules.}

We base, in part, our discussion of the general concept of
quantization and its application to homotopy of skein modules on
\cite{Tu-3,H-P-1}.
\begin{definition}\label{8.1}
\begin{enumerate}
\item [(a)] A Lie ring is a $Z$ module $B$, with a map
$[\ ,\ ] :B\times B \to B$ satisfying the condition:
\begin{enumerate}
\item [(i)] $[\ ,\ ]$ is $Z$ bilinear; that is:\\
$[x+y,z] =[x,z] + [y,z]$ and $[x,y+z]= [x,y]+[x,z]$.
\item [(ii)] $[x,x]=0$, in particular $[,]$ is anti-symmetric
($[x,y]=-[y,x]$).
\item [(iii)] $[\ ,\ ]$ satisfies the Jacobi identity:\\
$[x,[y,z]] + [y,[z,x]] + [z,[x,y]] =0$.
\end{enumerate}
\item [(b)] Let $R$ be a commutative ring and $R$-module $B$ 
a Lie ring, then $B$ is called an $R$-Lie algebra if\
$[ax,y]=[x,ay]=a[x,y]$, that is, $[,]$ is an $R$ bilinear map.
\end{enumerate}
\end{definition}
\begin{example}\label{8.2}
Let $B$ be a ring and define $[x,y]= xy-yx$. Then $B$ becomes
a Lie ring. If $B$ is an $R$-algebra then $(B,[,])$ becomes an
$R$-Lie algebra.
\end{example}
\begin{example}\label{8.3}
Consider the homotopy skein algebra ${\cal HS}(F\times I)$. As in
Example 8.2 it is a Lie algebra with a bracket defined by
$[L_{1},L_{2}]=L_{1}\cdot L_{2}-L_{2}\cdot L_{1}$. The Lie bracket
of knots is a linear combination of knots, thus the submodule 
of ${\cal HS}(F\times I)$
generated by knots, $Z[z]\hat\pi$, is a Lie subalgebra of ${\cal
HS}(F\times I)$. The bracket $[K_{1},K_{2}]$ can be written now as
$$[K_{1},K_{2}]=z\sum_{p\in K_{1}\cap K_{2}} sgn(p)(K_{1}\cdot
K_{2})^{p}_{0}.$$
As it is equal to $K_{1}\cdot K_{2}-K_{2}\cdot K_{1}$, the formula
does not depend on the choice of diagrams for $K_{1}$ and $K_{2}$.
The bracket was first considered by Goldman \cite{Gol} (for $z=1$). The
Lie algebra $Z\hat\pi$ is called the Goldman-Wolpert Lie
algebra of curves on $F$, and its relation to knot theory was
first noticed by Turaev \cite{Tu-2}.

\end{example}

\begin{example}\label{8.4}
Let $H$ be an abelian group (e.g. the first homology group of a surface)
with addition denoted by $\oplus$ , and consider a bilinear
anti-symmetric form\footnote{A form $f_{n}:H\times H\to Z_{n}$ and
the bracket $[\ ,\  ]_{q}:Z[q^{\pm 1}]/(\frac{q^n-q^{-n}}{q-q^{-1}})
H\times Z[q^{\pm 1}]/(\frac{q^n-q^{-n}}{q-q^{-1}})H
\to Z[q^{\pm 1}]/(\frac{q^n-q^{-n}}{q-q^{-1}})H$ can be considered
analogously.}
 on $H$, $f: H \times H \to Z$ (e.g. the homology intersection form). 
We have the structure of a Lie algebra on 
$Z[q^{\pm 1}]H$ with the bracket $[,]_{q}:Z[q^{\pm
1}]H\times Z[q^{\pm 1}]H\to Z[q^{\pm 1}]H$ defined on elements of
$H$ by $[g,h]_q = [f(g,h)]_q(g\oplus h)$, where the $q$-integer
$[n]_q$ is given by $[n]_q =\frac{q^n-q^{-n}}{q-q^{-1}}$, and
extended bilinearly to $Z[q^{\pm 1}]H$.\footnote{It may be more
``orthodox" to define a bracket by the formula
$[g,h]_{q}'=q^{f(g,h)}[f(g,h)]_{q}(g\oplus h)$ but this leads to a
$q$-Lie algebra. In particular
$[g,h]_{q}'=-q^{2f(g,h)}[h,g]_{q}'$.}
\\
 First notice that
$[-n]_q = -[n]_q$, so $[f(g,h)]_q = -[f(h,g)]_q$, and $[g,h]_q =
[f(g,h)]_q (g\oplus h) = -[f(h,g)]_q(h\oplus g)= -[h,g]_q$. We
will check the Jacobi identity on elements of $H$ using the
identity $[m+n]_{q}=q^{n}[m]_{q}+q^{-m}[n]_{q}.$\\ 
$[x,[y,z]_q]_q + [y,[z,x]_q]_q + [z,[x,y]_q]_q =$ \\
$[x,[f(y,z)]_q(y\oplus z)]_q + [y,[f(z,x)]_q(z\oplus x)]_q + 
[z,[f(x,y)]_q(x\oplus y)]_q =$\\
$[f(y,z)]_q[x,(y\oplus z)]_q + [f(z,x)]_q[y,(z\oplus x)]_q +
[f(x,y)]_q[z,(x\oplus y)]_q =$ \\
$[f(y,z)]_q[f(x,y\oplus z)]_q(x\oplus y\oplus z) + 
[f(z,x)]_q[f(y,z\oplus x)]_q (x\oplus y\oplus z) +$ \\
$ [f(x,y)]_q[f(z,x\oplus y)]_q (x\oplus y\oplus z)=
([f(y,z)]_q([f(x,y) + f(x,z)]_q)+$ \\
$[f(z,x)]_q([f(y,z)+f(y,x)]_q)+[f(x,y)]_q
([f(z,x)+f(z,y)]_q))(x\oplus y\oplus z)=$\\
$([f(y,z)]_q
(q^{f(x,z)}[f(x,y)]_q+q^{-f(x,y)}[f(x,z)]_q) +
[f(z,x)]_q(q^{f(y,x)}[f(y,z)]_q+ $ \\
$q^{-f(y,z)}[f(y,x)]_{q})+[f(x,y)]_{q}(q^{f(z,y)}
[f(z,x)]_{q}+q^{-f(z,x)}[f(z,y)]_{q}))(x\oplus y\oplus z)=$\\
$(q^{f(x,z)}[f(y,z)]_{q}[f(x,y)]_{q}+q^{-f(x,y)}[f(y,z)]_{q}[f(x,z)]_{q}
+q^{f(y,x)}[f(z,x)]_{q}[f(y,z)]_{q}+ $ \\
$q^{-f(y,z)}[f(z,x)]_{q}[f(y,x)]_{q}
+q^{f(z,y)}[f(x,y)]_{q}[f(z,x)]_{q}+q^{-f(z,x)}[f(x,y)]_{q}[f(z,y)]_{q})$ \\ 
$(x\oplus y\oplus z)=0.$
\end{example}

\begin{definition}\label{8.5}
A Poisson algebra is a commutative algebra equipped with a Lie
bracket which satisfies the following Leibniz rule:\\
$[ab,c]=a[b,c] + [a,c]b$.
\end{definition}

\begin{definition}\label{8.6}
\begin{enumerate}
\item[(a)] Let $B$ be an $R$-module and $B^{\otimes m}$ the tensor
product of $m$ copies of $B$ (with $B^{\otimes 0}=R$). Then the
tensor algebra ${\cal T}B$ is an $R$-module $\bigotimes_{i\geq 0}
B^{\otimes i}$ with the algebra multiplication defined by the
rule: $(a_1\otimes \ldots\otimes a_n)(b_1\otimes\ldots\otimes
b_m)=a_1\otimes\ldots\otimes a_n\otimes b_1\otimes\ldots\otimes
b_m$. If $B$ is a free module (with basis $E=\{e_i\}$) we can
identify ${\cal T}B$ with the algebra of noncommutative
polynomials in variables $\{e_i\}$, denoted by $R\{E\}$.
\item[(b)] The symmetric tensor algebra ${\cal S}B$ is the quotient
of ${\cal T}B$ by the ideal generated by commutators $a\otimes b-
b\otimes a$ where $a,b\in B$ (it suffices to consider a generating
set of $B$). If $B$ is a free module (with basis $E=\{e_i\}$) then
${\cal S}B$ is an algebra of (symmetric) polynomials in variables
$\{e_i\}$, $R[E]$.
\item[(c)] If $B$ is a Lie algebra then the universal enveloping algebra
${\cal U}B$ is the quotient of ${\cal T}B$ by the ideal generated
by expressions $a\otimes b- b\otimes a -[a,b]$ (${\cal U}B = {\cal
T}B/(a\otimes b- b\otimes a -[a,b])$).
\end{enumerate}
\end{definition}
\begin{example}
A symmetric tensor algebra ${\cal S}B$, with $B$ being an $R$-Lie
algebra is a Poisson algebra with the bracket extended from $B$ by
the Leibniz rule. We can write the global formula as follows:
$[a_1\otimes a_2\otimes \ldots \otimes a_m, b_1\otimes b_2\otimes
\ldots \otimes b_n] = \sum_{(i,j)} a_1\otimes a_2\otimes \ldots
\otimes a_{i-1}\otimes a_{i+1} \otimes \ldots \otimes a_m \otimes 
[a_i,b_j] \otimes b_1\otimes b_2\otimes \ldots \otimes
b_{j-1}\otimes b_{j+1} \otimes \ldots \otimes b_n$, where
$a_i,b_j\in B$.
\end{example}

Example 8.4 and topological motivation suggest the following
deformation of the Poisson algebra which we call a $q$-Poisson
algebra.

\begin{example}\label{8.8}
Let $Z[q^{\pm 1}]H$ with $[,]_{q}$ be the Lie algebra of Example
8.4. We define the $q$-symmetric tensor algebra as the quotient of
the tensor algebra by the ideal generated by relations
($q$-commutators) $q^{-f(g,h)}g\otimes h -q^{f(g,h)}h\otimes g$,
for $g,h\in H$. That is, ${\cal S}_{q}Z[q^{\pm 1}]H={\cal
T}Z[q^{\pm 1 }]H/(q^{-f(g,h)}g\otimes h -q^{f(g,h)}h\otimes g)$.

Let $a=g_{1}\otimes\ldots\otimes g_{k},$
$b=h_{1}\otimes\ldots\otimes h_{k'}$ and
$c=d_{1}\otimes\ldots\otimes d_{k''}$, where $g_{i},h_{i},d_{i}\in
H$. Define $f(a,b)=\sum_{i,j}f(g_{i},h_{j})$. We can extend our
bracket $[,]_{q}$ to ${\cal S}_{q}Z[q^{\pm 1}]H$ by a $q$-version
of the Leibniz rule:
$$[a\otimes b,c]_{q}=q^{-f(a,c)}a\otimes
[b,c]_{q}+q^{f(b,c)}[a,c]_{q}\otimes b.$$

Our Lie bracket $[,]_{q}$ on $Z[q^{\pm 1}]H$ leads also to a
''$q$-deformation" of an universal enveloping algebra:
$${\cal U}^{(q)}Z[q^{\pm 1}]H={\cal T}Z[q^{\pm
1}]H/(q^{-f(g,h)}g\otimes h -q^{f(g,h)}h\otimes g -[g,h]_{q}).$$

Our $q$-deformation of the universal enveloping algebra satisfies
the $q$-version of the Poincar{\'e}-Birkhoff-Witt theorem, that is,
${\cal U}^{(q)}Z[q^{\pm 1}]H$ is $Z[q^{\pm 1}]$-module isomorphic
to ${\cal S}_{q}Z[q^{\pm 1}]H.$\footnote{ Compute, in ${\cal
U}^{(q)}Z[q^{\pm 1}]H$, the balanced difference
$q^{-f(a,b)-f(a,c)-f(b,c)}a\otimes b\otimes
c-q^{f(a,b)+f(a,c)+f(b,c)}c\otimes b\otimes a$, for $a,b,c\in H$,
using two different methods (reflecting the symmetric group
relation $s_{1}s_{2}s_{1}=s_{2}s_{1}s_{2}$). The difference of
results:
$\Delta_{1}-\Delta_{2}=[[a,b]_q,c]_{q}+[[b,c]_q,a]_{q}+[[c,a]_q,b]_{q}$
is exactly the Jacobi expression which is equal to 0 as $Z[q^{\pm
1}]H$ is a Lie algebra (Example 8.4). This observation,
$\Delta_{1}-\Delta_{2}=0$, is the key in proving the
Poincar\'{e}-Birkhoff-Witt Theorem.}
\end{example}

 The relation between the
Poisson algebras defined above and homotopy skein modules (of
$F\times I$) is best formulated in the language of
``quantizations". In fact knot theory leads to several nontrivial
quantizations.

\begin{definition}\label{8.9}
\begin{enumerate}
\item [(a)]
Let $P$ be a Poisson algebra over $Z$ and let $A$ be an algebra
over $Z[z]$ which is free as a $Z[z]$-module. A $Z$-module
epimorphism $\phi: A \to P$ is called a Drinfeld-Turaev
quantization of $P$ if
\begin{enumerate}
\item [(i)]
$\phi(p(z)a)=p(0)\phi(a)$ for all $a\in A$ and all $p(z)\in Z[z]$,
and
\item [(ii)]
$ab-ba\in z{\phi}^{-1}([\phi(a),\phi(b)]) $ for all $a,b \in P$.
\end{enumerate}
\item [(b)] If we do not require, as in (a), that $A$ is free as
a $Z[z]$-module, we call this a {\it weak} Drinfeld-Turaev
quantization.
\end{enumerate}
\end{definition}
\begin{example}\label{8.10}
Let $A$ be an algebra over the polynomial ring $R[z]$ which is
free as an $R[z]$-module. Assume that the quotient algebra $A/zA$ 
is abelian so
that $ab-ba \in zA$. Then the formula $[a\ mod\ zA, b\ mod\ zA]
=z^{-1}(ab-ba)\ mod zA$ equips $A/zA$ with a Lie bracket which
satisfies the Leibniz rule. Thus $A \to A/zA$ is a quantization.
\end{example}

\begin{theorem}\label{8.11}
\begin{enumerate}
\item[(1.)] [\cite{H-P-1,Tu-3}]\ \\
Let $Q:{\cal HS}(F\times I)\to {\bf S}Z\hat\pi$ be the $Z$-module
homomorphism from the homotopy skein module of $F\times I$ to the
Poisson (symmetric tensor) algebra of the Lie algebra $Z\hat\pi$,
defined by $Q(\sum p_{i}(z)L_{i})=\sum p_{i}(0)Q(L_{i})$ where for
$L_{i}=K_{1}\cdot\ldots\cdot K_{m}\in{\cal B}(F)$, $Q(L_{i})=
=K_{1}\otimes\ldots\otimes K_{m}$. Then $Q$ is a quantization of
${\bf S}Z\hat\pi$.

\item[(2.)]
Let $V(M)$ be a submodule of ${\cal HS}^{q}(M)$ generated by
relations $K_{1}-K_{2}$, where $K_{1}$ and $K_{2}$ are homologous
knots in $M$ (that is, $K_{1}=K_{2}$ in $H_{1}(M)$). Then ${\cal
HS}^{q}(F\times I)/V(F\times I)$ is $Z[q^{\pm 1},z]$-algebra
isomorphic to ${\cal U}^{q,z}Z[q^{\pm 1},z]H_{1}(F\times I)={\cal
T}Z[q^{\pm 1},z]H_{1}(F\times I)/q^{-f(g,h)}g\otimes
h-q^{f(g,h)}h\otimes g-z[g,h]_{q})$, where $g,h\in H_{1}(F\times
I)$.

\item[(3.)]
If we relax the definition of quantization to allow deformation of
$q$-Poisson algebras (like in Example 8.8), then the map $Q:{\cal
U}^{q,z}Z[q^{\pm 1},z]H_{1}(F\times I)\to {\bf S}_{q}Z[q^{\pm
1}]H_{1}(F\times I)$ is a quantization.

\item[(4.)] If we modify the multiplication of links in $F \times I$ by
putting $L_1 \hat \cdot L_2 = q^{-f(L_1,L_2)}L_1 \cdot L_2$ where 
$f(L_1,L_2) = \Sigma_{p\in L_1\cap L_2} sgn(p)$, and the homomorphism
$Q$ so that $\hat Q$ is a $\hat\cdot$ algebra homomorphism 
from  ${\cal HS}^{q}(F\times I)/V(F\times I)$ 
to the symmetric tensor algebra ${\bf S}Z[q^{\pm 1}]H_{1}(F\times I)$, 
then $\hat Q$ is a quantization.
Furthermore, ${\cal HS}^{q}(F\times I)/V(F\times I)$ with the product 
$\hat \cdot$ is $Z[q^{\pm 1},z]$-algebra
isomorphic to ${\cal U}^{z}Z[q^{\pm 1},z]H_{1}(F\times I)={\cal
T}Z[q^{\pm 1},z]H_{1}(F\times I)/g\otimes h-h\otimes g-z[g,h]_{q})$, 
where $g,h\in H_{1}(F\times I)$.
 
\end{enumerate}
\end{theorem}

\begin{proof}
\begin{enumerate}
\item[(1.)]
It follows from Theorem 6.1(ii) that ${\cal HS}(F\times I)$ is
$Z[z]$-algebra isomorphic to ${\cal T}Z[z]\hat\pi /(K_{1}\otimes
K_{2}-K_{2}\otimes K_{1}-z[K_{1},K_{2}])$ where $[K_{1},K_{2}]$ is
the Goldman-Wolpert Lie bracket (Example 8.3). This algebra is in
turn a Drinfeld-Turaev quantization of the Poisson algebra ${\bf
S}Z\hat\pi$ (see Example 8.10).

\item[(2.)]
One should carefully follow the proof of Theorem 6.1. Notice that
(2.) reduces to Theorem 6.1(iii) for $\pi_{1}(F)$ abelian
($\pi_{1}(F)=H_{1}(F)$).

\item[(3.)]
For $z=0$,  ${\cal U}^{q,z}Z[q^{\pm 1},z]H_{1}(F\times I)$ reduces
to ${\bf S}_qZ[q^{\pm 1}]H_{1}(F\times I)$.

\item[(4.)] 
For $z=0$, ${\cal U}^{z}Z[q^{\pm 1},z]H_{1}(F\times I)$ reduces
to ${\bf S}Z[q^{\pm 1}]H_{1}(F\times I)$.   
The advantage of our modified product ($\hat \cdot$)
is that it allows us to equip the quotient 
${\cal HS}^{q}(F\times I)/V(F\times I)$
with the Hopf algebra structure 
(as ${\cal U}^{z}Z[q^{\pm 1},z]H_{1}(F\times I)$ is a Hopf algebra).
\end{enumerate}
\end{proof}

\section{Speculation}\label{9}
The only previous work on homotopy skein modules beyond
\cite{H-P-1,Tu-3} was the work by U.Kaiser
\cite{Kai}\footnote{There is related work by Andersen, Mattes and
Reshetikhin \cite{AMR-1,AMR-2} but it does not touch homotopy
skein modules.}. Kaiser showed that ${\cal HS}(M)$ is free for lens
spaces $L(p,q)$ ($p\neq 0$). One should be able to generalize his
result to the q-homotopy skein modules.

Of great interest is to produce invariants of 3-manifolds from
homotopy skein modules, say following the method of Lickorish
(compare also \cite{BHMV}). The first step is to extend homotopy skein
modules to framed links. The most natural solution is, following
Kauffman, to consider the skein relations of unframed links as
relations among framed links with zero framing (at least in a
homology spheres). Thus adding framing relation $L^{(1)}=sL$,
where $L^{(1)}$ denotes the framed link obtained from a framed
link $L$ by adding one positive twist to the framing of a
component of $L$. Our skein expressions then have the form:\\
$s^{-1}L_+^{fr} -sL_-^{fr}$ in the case of a self-crossing, and
$s^{-1}q^{-1}L_+^{fr} - sqL_-^{fr} -z L_0^{fr}$ in the case of a mixed
crossing. The resulting skein module is ${\cal HS}^{q,s}(M)= Z[q^{\pm
1},s^{\pm 1},z]{\cal L}^{fr}/ S_{q,s}$ where $S_{q,s}$ is a
submodule of $Z[q^{\pm 1},s^{\pm 1},z]{\cal L}^{fr}$ generated by
framing relations and homotopy skein relation as above. One can
hope that this skein module should produce several interesting
invariants of 3-manifolds including Dijkgraaf-Witten invariants
\cite{Di-Wi}.
One can get at least  the Murakami-Ohtsuki-Okada invariant \cite
{M-O-O}, as ${\cal HS}^{q,s}(M)$ dominates the second skein module
$S_{2}(M;q)$ which in turn can be used to construct the
Murakami-Ohtsuki-Okada invariant \cite{Pr-7}.

We can also work with homotopy skein modules of unoriented links.
In this case we can follow fruitfully the idea of D.Johnson
(\cite{Gol,Tu-2,H-P-1}) of embedding an unoriented link into the
sum of all its orientations. This can be used to analyze the
Kauffman homotopy skein modules (compare \cite{Tu-2,H-P-1}).
These modules, as well as related Vassiliev-Gusarov homotopy skein
modules (compare \cite{Pr-4}), are worthy of detailed
consideration.

\ \\ \ \\
\footnotesize{
\centerline{\it Department of Mathematics, University of Maryland}
\centerline{\it College Park, MD 20742}

The author is on leave from:\\
\centerline{\it Department of Mathematics, The George Washington University}
\centerline{\it 2201 G Str. Funger Hall, Washington, D.C. 20052}
\centerline{\it przytyck@gwu.edu}}
\end{document}